# SKEW CONVOLUTION SEMIGROUPS AND AFFINE MARKOV PROCESSES

BY D. A. DAWSON[1] AND ZENGHU LI[2]

*Carleton University and Beijing Normal University*

A general affine Markov semigroup is formulated as the convolution of a homogeneous one with a skew convolution semigroup. We provide some sufficient conditions for the regularities of the homogeneous affine semigroup and the skew convolution semigroup. The corresponding affine Markov process is constructed as the strong solution of a system of stochastic equations with non-Lipschitz coefficients and Poisson-type integrals over some random sets. Based on this characterization, it is proved that the affine process arises naturally in a limit theorem for the difference of a pair of reactant processes in a catalytic branching system with immigration.

**1. Introduction.** The concept of *affine processes* unifies a wide class of Markov processes including Ornstein–Uhlenbeck processes (OU-processes) and continuous state branching processes with immigration (CBI-processes). Those processes involve rich common mathematical structures and the unified treatment of them develops interesting connections between several areas in the theory of probability. The "affine property" is roughly that the logarithm of the characteristic function of the transition semigroup is given by an affine transformation of the initial state $x \mapsto \langle x, \psi(t,u) \rangle + \phi(t,u)$; see Section 3. An important special case is where the affine transformation is homogeneous, that is, it only contains a nontrivial linear part $\langle x, \psi(t,u) \rangle$. In this case, we refer to the affine semigroup as *homogeneous*. A general affine semigroup can be constructed as the convolution of a homogeneous one with an associated *skew convolution semigroup*, which corresponds to the constant term $\phi(t,u)$ and gives the one-dimensional distributions of the affine process

Received September 2004; revised May 2005.
[1]Supported by an NSERC Research grant.
[2]Supported by NSFC and NCET Research grants.
*AMS 2000 subject classifications.* Primary 60J35; secondary 60J80, 60H20, 60K37.
*Key words and phrases.* Skew convolution semigroup, affine process, continuous state branching process, catalytic branching process, immigration, Ornstein–Uhlenbeck process, stochastic integral equation, Poisson random measure.

started with the trivial initial state. A complete characterization of general finite-dimensional affine processes was recently given by Duffie, Filipović and Schachermayer [10] under a regularity assumption, which requires that the coefficients $\psi(t,u)$ and $\phi(t,u)$ are both differentiable at $t=0$. Based on this characterization, they discussed a wide range of applications of affine processes in mathematical finance.

The problems of characterizing different particular classes of affine processes have also been studied by some other authors. In particular, Watanabe [28] gave a complete description of regular two-dimensional continuous state branching processes. He also proved that the regularity property of such processes is implied by the stochastic continuity. A similar characterization of finite-dimensional continuous state branching processes was given in [23]. The same problem for measure-valued branching processes was investigated in [11]. In those cases, the processes are defined by homogeneous affine semigroups. On the other hand, a complete characterization for stochastically continuous one-dimensional CBI-processes was given in [16]; see also [27]. In the setting of measure-valued processes, Li [18] gave a formulation of immigration structures in terms of skew convolution semigroups. It was proved in [18] that the skew convolution semigroups associated with a measure-valued branching process are in 1–1 correspondence with a class of infinitely divisible probability entrance laws; see also [19, 20]. A construction of trajectories of the corresponding immigration processes was given in [21] by summing up measure-valued paths in some Kuznetsov processes. Skew convolution semigroups and OU-processes on real separable Hilbert spaces were studied in [3, 7, 8, 14, 22, 26]. Roughly speaking, a skew convolution semigroup is regular if and only if it is determined by a closable entrance law. For both the Hilbert spaces and the spaces of measures, a stochastically continuous skew convolution semigroup can be irregular. A number of such examples arising in applications were discussed in [7, 8, 19, 20].

The main purpose of this paper is to investigate the basic characterizations and regularities of affine Markov semigroups and processes. For simplicity of the presentation we shall confine ourselves to the nondegenerate two-dimensional case, but most of the arguments can be generalized to higher dimensions. From the results of [10] we know that the constant part in the affine structure is usually smoother than the linear part. Therefore, we discuss *separately* the regularities of homogeneous affine semigroups and those of skew convolution semigroups. It turns out that a skew convolution semigroup always consists of infinitely divisible probability measures. We prove that such a semigroup is regular if and only if the linear part of the logarithm of its characteristic function is absolutely continuous. Some sufficient conditions for the regularities of homogeneous affine semigroups and skew convolution semigroups are given in terms of their first moments. Those results give a partial solution of the problem of characterizing all affine semigroups without regularity assumption; see [10], Remark 2.11. We then give a construction of the affine process as the strong solution of a system of stochastic integral equations with random and non-Lipschitz coefficients and jumps of Poisson type selected from some random sets. A similar equation system is used to construct a class of catalytic CBI-processes. The concept of *catalytic branching processes* was first introduced by Dawson and Fleischmann [5] in the setting of measure-valued diffusions; see, for example, [4] and [6] for some related developments. As an application of the characterization by stochastic equations, we show that an affine process arises naturally in a limit theorem for the difference of a pair of reactant processes in a catalytic CBI-process. This result is of interest since it seems that the connection between affine processes and catalytic branching processes has not been noticed before. The studies of those two classes of processes have been undergoing rapid developments in recent years with different motivations. The interplay between them provides new motivations for both sides and might stimulate some further studies on related topics.

NOTATION. Let $\mathbb{R}_+ = [0, \infty)$ and $\mathbb{R}_- = (-\infty, 0]$. Let $\lambda$ denote the Lebesgue measure on $\mathbb{R}$. For a Borel measure $\nu$ and a Borel function $f$ on $E \subseteq \mathbb{R}$ or $\mathbb{R}^2$, we write $\nu(f)$ for $\int_E f \, d\nu$ if the integral exists. Write $\Delta_\xi f(x) = f(x + \xi) - f(x)$ if the right-hand side is well defined. Let $\hat{\nu}$ denote the *characteristic function* of $\nu$ defined by

$$\hat{\nu}(u) := \int_E \exp\{\langle u, \xi \rangle\} \nu(d\xi), \qquad u \in U,$$

where $U \subseteq \mathbb{C}$ or $\mathbb{C}^2$ is to be specified. For $x \in \mathbb{R}$ set $l_1(x) = |x|$ and $l_{12}(x) = |x| \wedge |x|^2$. Let

$$\chi(x) = \begin{cases} x, & \text{if } x \in [-1, 1], \\ 1, & \text{if } x \in (1, \infty), \\ -1, & \text{if } x \in (-\infty, 1). \end{cases}$$

For $x = (x_1, x_2) \in \mathbb{R}^2$ define $\chi(x) = (\chi(x_1), \chi(x_2))$. We make the convention that

$$\int_r^t = -\int_t^r = \int_{(r,t]} \quad \text{and} \quad \int_r^\infty = \int_{(r,\infty)}$$

for $r \leq t \in \mathbb{R}$.

**2. Homogeneous affine semigroups.** In this section we give the definition and prove some simple properties of homogeneous affine semigroups. Let $D = \mathbb{R}_+ \times \mathbb{R}$ and $U = \mathbb{C}_- \times (i\mathbb{R})$, where $\mathbb{C}_- = \{a + ib : a \in \mathbb{R}_-, b \in \mathbb{R}\}$ and $i\mathbb{R} = \{ib : b \in \mathbb{R}\}$. Note that the word "homogeneous" in the following definition has a meaning different from the one of "time-homogeneous."

DEFINITION 2.1. A transition semigroup $(Q(t))_{t\geq 0}$ with state space $D$ is called a *homogeneous affine semigroup* (*HA-semigroup*) if for each $t \geq 0$ there exists a continuous two-dimensional complex-valued function $\psi(t,\cdot) := (\psi_1(t,\cdot), \psi_2(t,\cdot))$ on $U$ such that

$$(2.1) \quad \int_D \exp\{\langle u, \xi \rangle\} Q(t, x, d\xi) = \exp\{\langle x, \psi(t, u)\rangle\}, \qquad x \in D, u \in U.$$

The HA-semigroup $(Q(t))_{t\geq 0}$ given by (2.1) is called *regular* if it is stochastically continuous and the derivative $\psi'_t(0, u)$ exists for all $u \in U$ and is continuous at $u = 0$.

PROPOSITION 2.1. *Let $(Q(t))_{t\geq 0}$ be a HA-semigroup defined by (2.1). Then $\psi(t, u) \in U$ and*

$$(2.2) \qquad \psi(r + t, u) = \psi(r, \psi(t, u)), \qquad r, t \geq 0, u \in U.$$

*Moreover, $\psi_2(t, u)$ has the form*

$$(2.3) \qquad \psi_2(t, u) = \beta_{22}(t) u_2, \qquad t \geq 0, u \in U,$$

*where $\beta_{22}(\cdot)$ is a function on $[0, \infty)$ satisfying*

$$(2.4) \qquad \beta_{22}(r + t) = \beta_{22}(r) \beta_{22}(t), \qquad r, t \geq 0.$$

PROOF. From (2.1) we see that $\langle x, \psi(t, u) \rangle \in \mathbb{C}_-$ for all $x \in D$. This implies that $\psi(t, u) \in U$ for all $u \in U$. For any $x_2 \in \mathbb{R}$ we have

$$Q(t, (0, x_2), \cdot) * Q(t, (0, -x_2), \cdot) = \delta_{(0,0)}.$$

Then $Q(t, (0, x_2), \cdot)$ must be degenerate. Let $Q(t, (0, x_2), \cdot) = \delta_{\beta(t, x_2)}$ for $\beta(t, x_2) \in D$. It follows that

$$(2.5) \quad \exp\{x_2 \psi_2(t, u)\} = \int_D \exp\{\langle u, \xi\rangle\} Q(t, (0, x_2), d\xi) = \exp\{\langle u, \beta(t, x_2)\rangle\}.$$

Then $x_2 \psi_2(t, u) = \langle u, \beta(t, x_2)\rangle$ and so $\beta(t, x_2) = \beta(t, 1) x_2$. Since $\beta(t, x_2) \in D$ for all $x_2 \in \mathbb{R}$, we must have $\beta_1(t, 1) = 0$ and hence $\psi_2(t, u) = \beta_2(t, 1) u_2$ for all $u \in U$. That is, (2.3) holds with $\beta_{22}(t) = \beta_2(t, 1)$. The relation (2.2) follows from (2.1) and the Chapman–Kolmogorov equation for $(Q(t))_{t\geq 0}$. By (2.2) and (2.3) we obtain

$$\beta_{22}(r + t) u_2 = \psi_2(r, \psi(t, u)) = \beta_{22}(r) \psi_2(t, u) = \beta_{22}(r) \beta_{22}(t) u_2,$$

which implies (2.4). □

COROLLARY 2.1. *If (2.1) defines a stochastically continuous HA-semigroup $(Q(t))_{t\geq 0}$, then there is a constant $\beta_{22} \in \mathbb{R}$ such that $\beta_{22}(t) = \exp\{\beta_{22} t\}$ for all $t \geq 0$.*

PROPOSITION 2.2. *Let $(Q(t))_{t\geq 0}$ be a HA-semigroup defined by* (2.1). *Then $\psi_1(t,u)$ has the representation*

(2.6)
$$\psi_1(t,u) = \beta_{11}(t)u_1 + \beta_{12}(t)u_2$$
$$+ \alpha(t)u_2^2 + \int_D (e^{\langle u,\xi\rangle} - 1 - u_2\chi(\xi_2))\mu(t,d\xi),$$

*where $\alpha(t) \in \mathbb{R}_+$, $(\beta_{11}(t), \beta_{12}(t)) \in D$ and $\mu(t, d\xi)$ is a $\sigma$-finite measure on $D$ supported by $D \setminus \{0\}$ such that*

$$\int_D [\chi(\xi_1) + \chi^2(\xi_2)]\mu(t,d\xi) < \infty.$$

*Moreover, for any $r, t \geq 0$ we have*

(2.7) $\quad\beta_{11}(r+t) = \beta_{11}(r)\beta_{11}(t),$

(2.8)
$$\beta_{12}(r+t) = \beta_{11}(r)\beta_{12}(t) + \beta_{12}(r)\beta_{22}(t)$$
$$+ \int_D [Q(t)\chi_2(\xi) - \beta_{22}(t)\chi(\xi_2)]\mu(r,d\xi),$$

(2.9) $\quad\alpha(r+t) = \beta_{11}(r)\alpha(t) + \alpha(r)\beta_{22}^2(t),$

(2.10) $\quad\mu(r+t,\cdot) = \int_D \mu(r,d\xi)Q(t,\xi,\cdot) + \beta_{11}(r)\mu(t,\cdot),$

*where*

(2.11) $\quad Q(t)\chi_2(\xi) = \int_D \chi(\eta_2)Q(t,\xi,d\eta)$

*and $\beta_{22}(t)$ is given by Proposition* 2.1.

PROOF. In view of (2.1) it is easy to see that each $Q(t,(x_1,0),\cdot)$ is an infinitely divisible probability measure on $D$. Then (2.6) follows by the special structure of $D$ and the Lévy–Khintchine representation for the characteristic function of an infinitely divisible distribution; see, for example, [17], pages 499–500. From (2.6) and the results of Proposition 2.1 we get

$$\psi_1(r+t,u) = \beta_{11}(r)\psi_1(t,u) + \beta_{12}(r)\beta_{22}(t)u_2 + \alpha(r)\beta_{22}^2(t)u_2^2$$
$$+ \int_D [e^{\langle \psi(t,u),\xi\rangle} - 1 - \beta_{22}(t)u_2\chi(\xi_2)]\mu(r,d\xi)$$
$$= \beta_{11}(r)\beta_{11}(t)u_1 + \beta_{11}(r)\beta_{12}(t)u_2$$
$$+ \beta_{11}(r)\alpha(t)u_2^2 + \beta_{12}(r)\beta_{22}(t)u_2$$
$$+ \alpha(r)\beta_{22}^2(t)u_2^2 + \beta_{11}(r)\int_D [e^{\langle u,\xi\rangle} - 1 - u_2\chi(\xi_2)]\mu(t,d\xi)$$
$$+ \int_D [e^{\langle \psi(t,u),\xi\rangle} - 1 - u_2 Q(t)\chi_2(\xi)]\mu(r,d\xi)$$

$$+ u_2 \int_D [Q(t)\chi_2(\xi) - \beta_{22}(t)\chi(\xi_2)]\mu(r, d\xi).$$

Then (2.7)–(2.10) follow by a comparison of the above expression with (2.6). □

LEMMA 2.1. *If* (2.1) *defines a stochastically continuous HA-semigroup and $\psi_1(t, u)$ is given by* (2.6), *then $t \mapsto \beta_{12}(t)$ is continuous.*

PROOF. The proof of [10], Lemma 3.1, shows that $\psi(t, u)$ is jointly continuous in $(t, u)$. Let $\mu_2(t, d\xi_2)$ denote the projection of $\mu(t, d\xi)$ to $\mathbb{R}$. Then

$$\psi_1(t, (0, iz)) = i\beta_{12}(t)z - \alpha(t)z^2 + \int_\mathbb{R} (e^{iz\xi_2} - 1 - iz\chi(\xi_2))\mu_2(t, d\xi_2)$$

is a continuous function of $(t, z) \in [0, \infty) \times \mathbb{R}$. It is not hard to find universal constants $0 < c_1 < c_2 < \infty$ so that

$$(2.12) \qquad c_1 \leq \left(1 - \frac{\sin \xi_2}{\xi_2}\right)\chi(\xi_2)^{-2} \leq c_2$$

for all $\xi_2 \in \mathbb{R} \setminus \{0\}$. Then for each $t \geq 0$ we can define the Borel measure $G(t, d\xi_2)$ on $\mathbb{R}$ by setting $G(t, \{0\}) = \alpha(t)/3$ and

$$G(t, d\xi_2) = \left(1 - \frac{\sin \xi_2}{\xi_2}\right)\mu_2(t, d\xi_2), \qquad \xi_2 \in \mathbb{R} \setminus \{0\}.$$

It follows that

$$(2.13) \quad \begin{aligned} \psi_1(t, (0, iz)) &= i\beta_{12}(t)z \\ &\quad + \int_\mathbb{R} (e^{iz\xi_2} - 1 - iz\chi(\xi_2))\left(1 - \frac{\sin \xi_2}{\xi_2}\right)^{-1} G(t, d\xi_2) \end{aligned}$$

where the integrand is defined at $\xi_2 = 0$ by continuity as $-3z^2$. By dominated convergence,

$$t \mapsto v(t, \lambda) := 2\psi_1(t, (0, i\lambda)) - \int_{\lambda-1}^{\lambda+1} \psi_1(t, (0, iz))\, dz$$

is continuous for each $\lambda \in \mathbb{R}$. One may check easily that $v(t, \lambda)$ is the characteristic function of $G(t, d\xi_2)$. Then Lévy's continuity theorem implies that $t \mapsto G(t, d\xi_2)$ is continuous by weak convergence. For any fixed $z \in \mathbb{R}$, the integrand in (2.13) is bounded and continuous in $\xi_2$, so the integral term is continuous in $t \geq 0$. By the continuity of $\psi_1(t, (0, iz))$ we find that $t \mapsto \beta_{12}(t)$ is continuous. □

PROPOSITION 2.3. *Suppose that* (2.1) *defines a stochastically continuous HA-semigroup and $\psi_1(t,u)$ is given by* (2.6). *Then*

$$(2.14) \quad B(t) := \beta_{11}(t) + \beta_{12}^2(t) + \alpha(t) + \int_D [\chi(\xi_1) + \chi^2(\xi_2)]\mu(t, d\xi)$$

*is a locally bounded function of $t \geq 0$.*

PROOF. Let $\mu_1(t, d\xi_1)$ denote the projection of $\mu(t, d\xi)$ to $\mathbb{R}_+$. Then

$$(2.15) \quad \psi_1(t, (-z, 0)) = -\beta_{11}(t)z + \int_0^\infty (e^{-z\xi_1} - 1)\mu_1(t, d\xi_1)$$

is continuous in $(t, z) \in [0, \infty) \times \mathbb{R}_+$. In particular, (2.15) is locally bounded in $t \geq 0$ for any fixed $z \in \mathbb{R}_+$. Taking $\lambda = 1$ one finds that

$$\beta_{11}(t) + \int_0^\infty \chi(\xi_1)\mu_1(t, d\xi_1)$$

is locally bounded in $t \geq 0$. By the proof of Lemma 2.1,

$$v(t, 0) = \frac{1}{3}\alpha(t) + \int_{\mathbb{R}} \left(1 - \frac{\sin \xi_2}{\xi_2}\right)\mu_2(t, d\xi_2)$$

is continuous and hence locally bounded in $t \geq 0$. In view of (2.12) we find that

$$(2.16) \quad \alpha(t) + \int_{\mathbb{R}} \chi^2(\xi_2)\mu_2(t, d\xi_2)$$

is also locally bounded in $t \geq 0$. By Lemma 2.1, $\beta_{12}(t)$ is locally bounded in $t \geq 0$. Then we have the desired result. $\square$

PROPOSITION 2.4. *Let $(Q(t))_{t \geq 0}$ be a stochastically continuous HA-semigroup defined by* (2.1). *Then there is a locally bounded nonnegative function $c_0(\cdot)$ on $[0, \infty)$ such that*

$$(2.17) \quad \int_D \chi(\xi_1)Q(t, x, d\xi) \leq c_0(t)\chi(x_1)$$

*and*

$$(2.18) \quad \int_D \chi^2(\xi_2)Q(t, x, d\xi) \leq c_0(t)[\chi(x_1) + \chi^2(x_2)].$$

PROOF. In view of (2.15), we have $\psi_1(t, (-z, 0)) \leq 0$. From (2.1) it follows that

$$\int_D (1 - e^{-z\xi_1})Q(t, x, d\xi) = 1 - \exp\{-x_1|\psi_1(t, (-z, 0))|\} \leq c_1(t, z)(1 - e^{-x_1}),$$

where $c_1(t,z) := 1 \vee |\psi_1(t,(-z,0))|$ is locally bounded in $t \geq 0$. Then we get the first inequality by letting $\lambda = 1$. The second inequality is obvious if $x_1 \geq 1$ or $|x_2| \geq 1$. On the other hand, we have

$$\int_D \left(1 - \frac{\sin \xi_2}{\xi_2}\right) Q(t,x,d\xi) = \frac{1}{2} \int_{-1}^1 dz \int_D (1 - e^{iz\xi_2}) Q(t,x,d\xi)$$

$$= \frac{1}{2} \int_{-1}^1 (1 - \exp\{x_1 \psi_1(t,(0,iz)) + ix_2 \beta_{22}(t)z\}) \, dz$$

$$\leq \frac{1}{2} \int_{-1}^1 [x_1 |\psi_1(t,(0,iz))| + h(t,x_1,x_2,z)] \, dz,$$

where

$$h(t,x_1,x_2,z) = \sum_{k=2}^\infty \frac{1}{k!} (x_1 |\psi_1(t,(0,iz))| + |x_2 \beta_{22}(t)z|)^k$$

$$\leq \sum_{k=2}^\infty \frac{1}{k!} (x_1 + |x_2|)^k (|\psi_1(t,(0,iz))| + |\beta_{22}(t)z|)^k$$

$$\leq 2(x_1^2 + x_2^2) q(t,x_1,x_2,z)$$

with

$$q(t,x_1,x_2,z) = \sum_{k=2}^\infty \frac{1}{k!} (x_1 + |x_2|)^{k-2} (|\psi_1(t,(0,iz))| + |\beta_{22}(t)z|)^k.$$

It follows that

$$\int_D \left(1 - \frac{\sin \xi_2}{\xi_2}\right) Q(t,x,d\xi) \leq \frac{1}{2} x_1 \int_{-1}^1 |\psi_1(t,(0,iz))| \, dz$$

$$+ (x_1^2 + x_2^2) \int_{-1}^1 q(t,x_1,x_2,z) \, dz,$$

which implies (2.18) for $0 \leq x_1 \leq 1$ and $|x_2| \leq 1$. □

As an immediate consequence of the above proposition we have the following.

COROLLARY 2.2. *Suppose that $(Q(t))_{t\geq 0}$ is a stochastically continuous HA-semigroup defined by (2.1). Let $U_\varepsilon = [0,\varepsilon) \times (-\varepsilon,\varepsilon)$ for $\varepsilon > 0$. Then for each $T \geq 0$ we have*

(2.19) $$\lim_{|x|\to 0} \sup_{0\leq t\leq T} Q(t,x,D\setminus U_\varepsilon) = 0.$$

**3. Skew convolution semigroups.** In this section we give a formulation of the general affine Markov semigroup in terms of a homogeneous one and a skew convolution semigroup. It turns out that a skew convolution semigroup always consists of infinitely divisible probability measures. We prove that such a semigroup is regular if and only if the linear part of the logarithm of its characteristic function is absolutely continuous. We shall fix a regular HA-semigroup $(Q(t))_{t\geq 0}$ on $D$ defined by (2.1), where $\psi(t,u) = (\psi_1(t,u), \beta_{22}(t)u_2)$ is given by Corollary 2.1 and Proposition 2.2.

DEFINITION 3.1. A family of probability measures $(\gamma(t))_{t\geq 0}$ on $D$ is called a *skew convolution semigroup* (*SC-semigroup*) associated with $(Q(t))_{t\geq 0}$ if

(3.1) $$\gamma(r+t) = (\gamma(r)Q(t)) * \gamma(t), \qquad r, t \geq 0,$$

where "$*$" denotes the convolution operation and $\gamma(r)Q(t)$ is the probability measure on $D$ defined by

(3.2) $$\gamma(r)Q(t)(B) = \int_D Q(t,\xi,B)\gamma(r,d\xi), \qquad B \in \mathscr{B}(D).$$

The concept of SC-semigroup generalizes that of the usual convolution semigroup; see also [7, 18, 19, 20, 21]. We refer the reader to [2, 24] for the general theory of convolution semigroups and Lévy processes.

PROPOSITION 3.1. *Let $(\gamma(t))_{t\geq 0}$ be a stochastically continuous SC-semigroup associated with $(Q(t))_{t\geq 0}$. Then each $\gamma(t)$ is an infinitely divisible probability measure, so we have the representation*

(3.3) $$\int_D \exp\{\langle u, \xi\rangle\}\gamma(t,d\xi) = \exp\{\phi(t,u)\}, \qquad u \in U,$$

*with*

(3.4) $$\phi(t,u) = b_1(t)u_1 + b_2(t)u_2 + a(t)u_2^2$$
$$+ \int_D (e^{\langle u,\xi\rangle} - 1 - u_2\chi(\xi_2))m(t,d\xi),$$

*where $a(t) \in \mathbb{R}_+$, $(b_1(t), b_2(t)) \in D$ and $m(t, d\xi)$ is a $\sigma$-finite measure on $D$ supported by $D \setminus \{0\}$ such that*

$$\int_D [\chi(\xi_1) + \chi^2(\xi_2)]m(t,d\xi) < \infty.$$

*Moreover, for any $r, t \geq 0$ we have*

(3.5) $$b_1(r+t) = b_1(r)\beta_{11}(t) + b_1(t),$$
$$b_2(r+t) = b_1(r)\beta_{12}(t) + b_2(r)\beta_{22}(t) + b_2(t)$$

(3.6)
$$+ \int_D [Q(t)\chi_2(\xi) - \beta_{22}(t)\chi(\xi_2)] m(r, d\xi),$$

(3.7) $\qquad a(r+t) = b_1(r)\alpha(t) + a(r)\beta_{22}^2(t) + a(t),$

(3.8) $\qquad m(r+t, \cdot) = \int_D m(r, d\xi) Q(t, \xi, \cdot) + b_1(r)\mu(t, \cdot) + m(t, \cdot).$

PROOF. Based on Corollary 2.2, the proof is a simplification of the arguments of Schmuland and Sun [26]. Let $t \geq 0$ be fixed. For each integer $n \geq 1$ we may use (3.1) inductively to obtain

$$\gamma(t) = \prod_{j=1}^{n} *\gamma(t/n) Q((j-1)t/n).$$

From the stochastic continuity of the SC-semigroup, we have $\lim_{n \to 0} \gamma(t/n) = \delta_0$. By virtue of Corollary 2.2, it is easy to show that $\{\gamma(t/n)Q((j-1)t/n) : j = 1, \ldots, n; n = 1, 2, \ldots\}$ form an infinitesimal triangular array. It follows that $\gamma(t)$ is infinitely divisible. Then we have representation (3.3) with $\phi(t, u)$ given by (3.4); see, for example, [17], pages 499–500, 515–519. From (3.1) we have

(3.9) $\qquad \phi(r+t, u) = \phi(r, \psi(t, u)) + \phi(t, u), \qquad r, t \geq 0, u \in U.$

Then relations (3.5)–(3.8) follow as in the proof of Proposition 2.2. $\square$

DEFINITION 3.2. A SC-semigroup $(\gamma(t))_{t \geq 0}$ associated with $(Q(t))_{t \geq 0}$ is called *regular* if $\phi(t, u) = \log \hat{\gamma}(t, u)$ has representation

(3.10) $\qquad \phi(t, u) = \int_0^t F(\psi(s, u)) \, ds, \qquad t \geq 0, u \in U,$

where $F = \log \hat{\nu}$ for an infinitely divisible probability measure $\nu$ on $D$.

We remark that if $\nu$ is an infinitely divisible probability measure on $D$, the function $F$ is well defined and (3.10) really determines a SC-semigroup. A simple but irregular SC-semigroup can be constructed by letting $Q(t)$ be the identity and letting $\gamma(t) = \delta_{(0, b_2(t))}$, where $b_2(t)$ is a discontinuous solution of $b_2(r+t) = b_2(r) + b_2(t)$; see, for example, [24], page 37. This example shows that some condition on the function $t \mapsto b_2(t)$ has to be imposed to get the regularity of the SC-semigroup $(\gamma(t))_{t \geq 0}$ given by (3.3) and (3.4).

DEFINITION 3.3. A transition semigroup $(P(t))_{t \geq 0}$ on $D$ is called a (general) *affine semigroup* associated with the HA-semigroup $(Q(t))_{t \geq 0}$ if its characteristic function has representation

(3.11) $\qquad \int_D \exp\{\langle u, \xi \rangle\} P(t, x, d\xi) = \exp\{\langle x, \psi(t, u) \rangle + \phi(t, u)\},$

where $\phi(t,\cdot)$ is a continuous function on $U$ satisfying $\phi(t,0) = 0$. The affine semigroup $(P(t))_{t\geq 0}$ defined by (3.11) is called *regular* if it is stochastically continuous and the derivatives $\psi'_t(0,u)$ and $\phi'_t(0,u)$ exist for all $u \in U$ and are continuous at $u = 0$.

PROPOSITION 3.2. *Let $(\gamma(t))_{t\geq 0}$ be a stochastically continuous SC-semigroup associated with the HA-semigroup $(Q(t))_{t\geq 0}$. Then $P(t,x,\cdot) = Q(t,x,\cdot) * \gamma(t,\cdot)$ defines a Feller affine semigroup $(P(t))_{t\geq 0}$.*

PROOF. It is easy to show that the kernels $P(t,x,\cdot)$ satisfy the Chapman–Kolmogorov equation. From [10], Proposition 8.2 we know that $(Q(t))_{t\geq 0}$ is a Feller semigroup. For any $f \in C_0(D)$ we have

$$P(t)f(x) = \int_D \gamma(t,dy) \int_D f(\xi+y) Q(t,x,d\xi).$$

Then we can use dominated convergence to find that $P(t)f \in C_0(D)$. Since both $(Q(t))_{t\geq 0}$ and $(\gamma(t))_{t\geq 0}$ are stochastically continuous, so is $(P(t))_{t\geq 0}$. It follows that $(P(t))_{t\geq 0}$ is a Feller semigroup. □

Clearly, if $(\gamma(t))_{t\geq 0}$ is a regular SC-semigroup, then the general affine semigroup $(P(t))_{t\geq 0}$ defined in Proposition 3.2 is also regular. To study the regularity of SC-semigroups, we need some preliminary results. The proofs in the sequel rely heavily on estimates derived from the relations (3.5)–(3.8).

PROPOSITION 3.3. *Suppose that (3.3) and (3.4) define a stochastically continuous SC-semigroup $(\gamma(t))_{t\geq 0}$. Then*

$$(3.12) \quad A(t) := b_1(t) + b_2^2(t) + a(t) + \int_D [\chi(\xi_1) + \chi^2(\xi_2)] m(t,d\xi)$$

*is a locally bounded function of $t \geq 0$. Moreover, we have $A(t) \to 0$ as $t \to 0$.*

PROOF. The stochastic continuity of the SC-semigroup implies that $\phi(t,u)$ is jointly continuous in $(t,u)$. Since $\phi(t,u) \to 0$ as $t \to 0$, the results follow by slight modifications of the arguments in the proof of Proposition 2.3. □

Let $B(\cdot)$ and $c_0(\cdot)$ be given respectively by Propositions 2.3 and 2.4. In the next two lemmas, we fix a constant $T \geq 0$ and let $C(T) \geq 0$ be a constant such that

$$(3.13) \quad \max\{B(t), c_0(t), \beta_{22}^2(t)\} \leq C(T), \qquad 0 \leq t \leq T.$$

For $0 \leq r_1 < t_1 < r_2 < t_2 < \cdots$ we set $\sigma_n = \sum_{j=1}^n (t_j - r_j)$.

LEMMA 3.1. *Suppose that* (3.3) *and* (3.4) *define a stochastically continuous SC-semigroup* $(\gamma(t))_{t\geq 0}$. *Then for* $0 \leq r_1 < t_1 < r_2 < t_2 < \cdots \leq T$ *we have*

$$\sum_{j=1}^{n}[b_1(t_j) - b_1(r_j)] \leq C(T)b_1(\sigma_n) \tag{3.14}$$

*and*

$$\sum_{j=1}^{n}[a(t_j) - a(r_j)] \leq C(T)[b_1(\sigma_n) + a(\sigma_n)]. \tag{3.15}$$

*Consequently, $b_1(t)$ and $a(t)$ are absolutely continuous in $t \geq 0$.*

PROOF. We shall only give the proof of (3.15) since the proof of (3.14) is similar. By (3.7) we find that $t \mapsto a(t)$ is nondecreasing and

$$a(t_1) - a(r_1) = b_1(t_1 - r_1)\alpha(r_1) + a(t_1 - r_1)\beta_{22}^2(r_1)$$
$$\leq C(T)[b_1(t_1 - r_1) + a(t_1 - r_1)],$$

that is, (3.15) holds for $n = 1$. Now suppose that (3.15) is true for $n-1$. By (3.13) and Propositions 2.2 and 3.1,

$$\sum_{j=1}^{n}[a(t_j) - a(r_j)] \leq [a(t_n) - a(r_n)] + C(T)[b_1(\sigma_{n-1}) + a(\sigma_{n-1})]$$

$$= b_1(t_n - r_n)\alpha(r_n) + a(t_n - r_n)\beta_{22}^2(r_n)$$
$$\quad + C(T)[b_1(\sigma_{n-1}) + a(\sigma_{n-1})]$$
$$= b_1(t_n - r_n)[\beta_{11}(\sigma_{n-1})\alpha(r_n - \sigma_{n-1})$$
$$\quad + \alpha(\sigma_{n-1})\beta_{22}^2(r_n - \sigma_{n-1})]$$
$$\quad + a(t_n - r_n)\beta_{22}^2(r_n) + C(T)[b_1(\sigma_{n-1}) + a(\sigma_{n-1})]$$
$$\leq C(T)b_1(t_n - r_n)\beta_{11}(\sigma_{n-1}) + C(T)b_1(t_n - r_n)\alpha(\sigma_{n-1})$$
$$\quad + C(T)[a(t_n - r_n)\beta_{22}^2(\sigma_{n-1}) + a(\sigma_{n-1})] + C(T)b_1(\sigma_{n-1})$$
$$= C(T)b_1(\sigma_n) + C(T)a(\sigma_n).$$

That proves (3.15) by induction. The absolute continuity of $b_1(t)$ and $a(t)$ follows by Proposition 3.3. □

LEMMA 3.2. *Suppose that* (3.3) *and* (3.4) *define a stochastically continuous SC-semigroup* $(\gamma(t))_{t\geq 0}$. *Set*

$$f(t) = \int_D \chi(\xi_1) m(t, d\xi) \quad \text{and} \quad g(t) = \int_D \chi^2(\xi_2) m(t, d\xi). \tag{3.16}$$

*Then for any $0 \leq r_1 < t_1 < \cdots < r_n < t_n \leq T$ we have*

(3.17) $$\sum_{j=1}^{n}[f(t_j) - f(r_j)] \leq C(T)[b_1(\sigma_n) + f(\sigma_n)]$$

*and*

(3.18) $$\sum_{j=1}^{n}[g(t_j) - g(r_j)] \leq C(T)[b_1(\sigma_n) + f(\sigma_n) + g(\sigma_n)].$$

*Consequently, $f(t)$ and $g(t)$ are absolutely continuous in $t \geq 0$.*

PROOF. The proofs of (3.17) and (3.18) are based on ideas similar to those in the proof of the last lemma. We here give the proof of (3.18) since it involves more careful calculations. By (3.8) we find that $t \mapsto g(t)$ is nondecreasing and

$$\begin{aligned}
g(t_1) - g(r_1) &= \int_D m(t_1 - r_1, d\xi) \int_D \chi^2(\eta_2) Q(r_1, \xi, d\eta) \\
&\quad + b_1(t_1 - r_1) \int_D \chi^2(\xi_2) \mu(r_1, d\xi) \\
&\leq C(T) \int_D [\chi(\xi_1) + \chi^2(\xi_2)] m(t_1 - r_1, d\xi) + C(T) b_1(t_1 - r_1) \\
&= C(T)[b_1(t_1 - r_1) + f(t_1 - r_1) + g(t_1 - r_1)],
\end{aligned}$$

where we used Proposition 2.4 for the inequality. Then (3.18) holds for $n = 1$. Suppose the inequality is true for $n - 1$. By (3.13) and the results of Propositions 2.4 and 3.1 we have

$$\begin{aligned}
\sum_{j=1}^{n}[g(t_j) - g(r_j)] &\leq [g(t_n) - g(r_n)] + C(T)[b_1(\sigma_{n-1}) + f(\sigma_{n-1}) + g(\sigma_{n-1})] \\
&= \int_D m(t_n - r_n, d\xi) \int_D \chi^2(\eta_2) Q(r_n, \xi, d\eta) \\
&\quad + b_1(t_n - r_n) \int_D \chi^2(\xi_2) \mu(r_n, d\xi) \\
&\quad + C(T)[b_1(\sigma_{n-1}) + f(\sigma_{n-1}) + g(\sigma_{n-1})] \\
&\leq C(T) \int_D m(t_n - r_n, d\xi) \int_D [\chi(\eta_1) + \chi^2(\eta_2)] Q(\sigma_{n-1}, \xi, d\eta) \\
&\quad + b_1(t_n - r_n) \int_D \mu(\sigma_{n-1}, d\xi) \\
&\quad \times \int_D \chi^2(\eta_2) Q(r_n - \sigma_{n-1}, \xi, d\eta)
\end{aligned}$$

$$+ b_1(t_n - r_n)\beta_{11}(\sigma_{n-1}) \int_D \chi^2(\eta_2)\mu(r_n - \sigma_{n-1}, d\eta)$$
$$+ C(T)[b_1(\sigma_{n-1}) + f(\sigma_{n-1}) + g(\sigma_{n-1})]$$
$$\leq C(T) \int_D m(t_n - r_n, d\xi) \int_D [\chi(\eta_1) + \chi^2(\eta_2)] Q(\sigma_{n-1}, \xi, d\eta)$$
$$+ C(T)b_1(t_n - r_n) \int_D [\chi(\xi_1) + \chi^2(\xi_2)]\mu(\sigma_{n-1}, d\xi)$$
$$+ C(T)[b_1(\sigma_n) + f(\sigma_{n-1}) + g(\sigma_{n-1})]$$
$$= C(T)[b_1(\sigma_n) + f(\sigma_n) + g(\sigma_n)].$$

Then (3.18) follows by induction. The second assertion follows by Proposition 3.3. □

LEMMA 3.3. *Suppose that* (3.3) *and* (3.4) *determine a stochastically continuous SC-semigroup* $(\gamma(t))_{t\geq 0}$. *Then there is a $\sigma$-finite kernel* $m'(s, d\xi)$ *from* $(0, \infty)$ *to* $D$ *supported by* $D \setminus \{0\}$ *so that*

$$(3.19) \qquad m(t, d\xi) = \int_0^t m'(s, d\xi) \, ds, \qquad t \geq 0.$$

PROOF. From (3.8) we see that $m(t, d\xi)$ is increasing in $t \geq 0$. By Lemma 3.2, the function

$$t \mapsto \int_D [\chi(\xi_1) + \chi^2(\xi_2)] m(t, d\xi)$$

is absolutely continuous. Then the assertion follows as in the proof of [8], Theorem 2.2. □

THEOREM 3.1. *Let* $(\gamma(t))_{t\geq 0}$ *be a stochastically continuous SC-semigroup given by* (3.3) *and* (3.4). *Then* $(\gamma(t))_{t\geq 0}$ *is regular if and only if the function* $t \mapsto b_2(t)$ *is absolutely continuous on* $[0, \infty)$.

PROOF. Suppose that $t \mapsto b_2(t)$ is absolutely continuous. In view of Lemma 3.1, we can find Borel measurable functions $a'(\cdot) \geq 0$ and $b'_j(\cdot)$ such that

$$a(t) = \int_0^t a'(s) \, ds \quad \text{and} \quad b_j(t) = \int_0^t b'_j(s) \, ds, \qquad j = 1, 2.$$

Let $m'(s, d\xi)$ be given by Lemma 3.3 and let

$$\phi'(s, u) = b'_1(s)u_1 + b'_2(s)u_2 + a'(s)u_2^2 + \int_D (e^{\langle u, \xi \rangle} - 1 - \chi(\xi_2)u_2) m'(s, d\xi).$$

Then we have

$$(3.20) \qquad \phi(t, u) = \int_0^t \phi'(s, u) \, ds.$$

By (3.9) and (3.20) it is easy to show that

$$\int_0^r \phi'(s+t,u)\,ds = \int_0^r \phi'(s,\psi(t,u))\,ds. \tag{3.21}$$

Let $\nu_s$ be the infinitely divisible probability measure on $D$ such that $\hat{\nu}_s(u) = \exp\{\phi'(s,u)\}$. Based on (3.21), it is easy to modify the definitions of $\phi'(s,\cdot)$ and $\nu_s$ accordingly so that $\nu_t = \nu_s Q_{t-s}$ for all $t > s > 0$ while (3.20) remains true; see [18]. In other words, $(\nu_s)_{s>0}$ is an entrance law for $(Q(t))_{t\geq 0}$. But $(Q(t))_{t\geq 0}$ is a Feller semigroup by [10], Proposition 8.2. Then the Ray–Knight compactification of $D$ with respect to $(Q(t))_{t\geq 0}$ coincides with its one-point compactification $\bar{D} := D \cup \{\partial\}$ and the Ray–Knight extension of $(Q(t))_{t\geq 0}$ satisfies $Q(t,\partial,\cdot) = \delta_\partial$ and $Q(t,x,\{\partial\}) = 0$ for every $x \in D$. It follows that the entrance space for $(Q(t))_{t\geq 0}$ is just $D$. By [25], page 196, there is a probability measure $\nu_0$ on $D$ such that $\nu_s = \nu_0 Q_s$ for all $s > 0$. Then $\phi'(s,u) = \log \hat{\nu}_0(\psi(s,u))$ and hence $(\gamma(t))_{t\geq 0}$ is regular. Conversely, if $(\gamma(t))_{t\geq 0}$ is regular, the function $t \mapsto \phi(t,u)$ is absolutely continuous on $[0,\infty)$ for every $u \in U$. Then (3.4) and Lemmas 3.1 and 3.3 imply that $t \mapsto b_2(t)$ is absolutely continuous. $\square$

COROLLARY 3.1. *Let $(\gamma(t))_{t\geq 0}$ be a stochastically continuous SC-semigroup defined by* (3.1). *Then $(\gamma(t))_{t\geq 0}$ is regular if and only if $t \mapsto \hat{\gamma}(t,u)$ is absolutely continuous for all $u \in U$.*

PROOF. By (3.4) and Lemmas 3.1 and 3.3, $t \mapsto b_2(t)$ is absolutely continuous if and only if $t \mapsto \phi(t,u)$ is absolutely continuous on $[0,\infty)$ for all $u \in U$. Then the result follows from Theorem 3.1. $\square$

**4. Regularities under moment conditions.** In this section we prove the regularities of HA-semigroups and their associated SC-semigroups under some conditions on the first moments. Suppose that $(Q(t))_{t\geq 0}$ is a stochastically continuous HA-semigroup defined by (2.1), where $\psi(t,u) = (\psi_1(t,u), \beta_{22}(t)u_2)$ is given by Corollary 2.1 and Proposition 2.2. Let us consider the following hypothesis.

HYPOTHESIS 4.1. Suppose that

$$\int_D [\xi_1 + l_{12}(\xi_2)]\mu(t,d\xi) < \infty \tag{4.1}$$

for all $t \geq 0$ or, equivalently,

$$\int_D |\xi| Q(t,x,d\xi) < \infty \tag{4.2}$$

for all $t \geq 0$ and $x \in D$.

If Hypothesis 4.1 holds, we have a more convenient representation for the function $\psi_1(t, u)$. Indeed, we may differentiate both sides of (2.6) to see that

$$(4.3) \qquad q_{11}(t) := \frac{\partial \psi_1}{\partial u_1}(t, u)\bigg|_{u=0} = \beta_{11}(t) + \int_D \xi_1 \mu(t, d\xi)$$

and

$$(4.4) \qquad q_{12}(t) := \frac{\partial \psi_1}{\partial u_2}(t, u)\bigg|_{u=0} = \beta_{12}(t) + \int_D [\xi_2 - \chi(\xi_2)] \mu(t, d\xi).$$

On the other hand, differentiating both sides of (2.1) we find that

$$(4.5) \qquad \int_D \xi_1 Q(t, x, d\xi) = x_1 q_{11}(t)$$

and

$$(4.6) \qquad \int_D \xi_2 Q(t, x, d\xi) = x_1 q_{12}(t) + x_2 \beta_{22}(t).$$

PROPOSITION 4.1. *If Hypothesis 4.1 holds, we have*

$$(4.7) \qquad \begin{aligned} \psi_1(t, u) &= \beta_{11}(t) u_1 + q_{12}(t) u_2 \\ &\quad + \alpha(t) u_2^2 + \int_D (e^{\langle u, \xi \rangle} - 1 - u_2 \xi_2) \mu(t, d\xi), \end{aligned}$$

*where $\alpha(t)$, $\beta_{11}(t)$ and $\mu(t, d\xi)$ are as in Proposition 3.1 and $q_{12}(t)$ is given by (4.4) and satisfies*

$$(4.8) \qquad q_{12}(r + t) = q_{11}(r) q_{12}(t) + q_{12}(r) \beta_{22}(t), \qquad r, t \geq 0.$$

PROOF. The representation (4.7) follows immediately from (2.6). By Proposition 2.2,

$$\begin{aligned}
q_{12}(r + t) &= \beta_{12}(r + t) + \int_D [\xi_2 - \chi(\xi_2)] \mu(r + t, d\xi) \\
&= \beta_{11}(r) \beta_{12}(t) + \beta_{12}(r) \beta_{22}(t) \\
&\quad + \int_D [Q(t) \chi_2(\xi) - \beta_{22}(t) \chi(\xi_2)] \mu(r, d\xi) \\
&\quad + \int_D [\xi_1 q_{12}(t) + \xi_2 \beta_{22}(t) - Q(t) \chi_2(\xi)] \mu(r, d\xi) \\
&\quad + \beta_{11}(r) \int_D [\xi_2 - \chi(\xi_2)] \mu(t, d\xi) \\
&= \beta_{11}(r) \beta_{12}(t) + \beta_{12}(r) \beta_{22}(t) + \beta_{22}(t) \int_D [\xi_2 - \chi(\xi_2)] \mu(r, d\xi) \\
&\quad + q_{12}(t) \int_D \xi_1 \mu(r, d\xi) + \beta_{11}(r) \int_D [\xi_2 - \chi(\xi_2)] \mu(t, d\xi) \\
&= \beta_{11}(r) q_{12}(t) + q_{12}(r) \beta_{22}(t) + q_{12}(t) \int_D \xi_1 \mu(r, d\xi).
\end{aligned}$$

Then we get (4.8) from (4.3). □

HYPOTHESIS 4.2. *Suppose that* (4.2) *holds and for every fixed* $x \in D$ *the mapping*

$$(4.9) \qquad t \mapsto |\xi| Q(t, x, d\xi)$$

*is continuous by the weak convergence of finite measures.*

The following theorem can be regarded as an extension of Watanabe ([28], Theorem 5) to the state space of the positive half plane.

THEOREM 4.1. *If Hypothesis* 4.2 *holds, the HA-semigroup* $(Q(t))_{t \geq 0}$ *is regular.*

PROOF. Under the hypothesis, we may differentiate both sides of (2.1) with respect to $u_1$ and $u_2$ to get

$$\int_D \xi_j \exp\{\langle u, \xi \rangle\} Q(t, (x_1, 0), d\xi) = x_1 \psi'_{1,u_j}(t, u) \exp\{x_1 \psi_1(t, u)\}$$

for $j = 1, 2$. On the other hand, since (4.9) depends on $t \geq 0$ continuously, we have

$$\lim_{t \to 0} \int_D \xi_1 \exp\{\langle u, \xi \rangle\} Q(t, (x_1, 0), d\xi) = x_1 \exp\{u_1 x_1\}$$

and

$$\lim_{t \to 0} \int_D \xi_2 \exp\{\langle u, \xi \rangle\} Q(t, (x_1, 0), d\xi) = 0.$$

Comparing the above equalities we obtain

$$(4.10) \qquad \lim_{t \to 0} \psi'_{1,u_1}(t, u) = 1 \quad \text{and} \quad \lim_{t \to 0} \psi'_{1,u_2}(t, u) = 0.$$

For $t > 0$ let

$$p(t, u) = \frac{1}{t} \int_0^t \psi(s, u) \, ds.$$

Then (4.10) implies that $\lim_{t \to 0} p'_{1,u_1}(t, u) = 1$ and $\lim_{t \to 0} p'_{1,u_2}(t, u) = 0$. By Proposition 2.1 and Corollary 2.1 we have $p'_{2,u_1}(t, u) = 0$ and $\lim_{t \to 0} p'_{2,u_2}(t, u) = 1$. Let $U_1$ be any fixed bounded subset of $U$. It is easy to see that the above limits hold with uniform convergence on $U_1$. Then we can choose sufficiently small $r > 0$ so that the matrix

$$(4.11) \qquad \partial p(r, u) := \begin{pmatrix} p'_{1,u_1}(r, u) & p'_{2,u_1}(r, u) \\ p'_{1,u_2}(r, u) & p'_{2,u_2}(r, u) \end{pmatrix}$$

is invertible for all $u \in U_1$. Observe that

$$p(r, \psi(t,u)) - p(r,u) = \frac{1}{r}\left[\int_0^r \psi(s, \psi(t,u))\,ds - \int_0^r \psi(s,u)\,ds\right]$$
$$= \frac{1}{r}\left[\int_t^{r+t} \psi(s,u)\,ds - \int_0^r \psi(s,u)\,ds\right]$$
$$= \frac{1}{r}\int_0^t [\psi(r+s, u) - \psi(s,u)]\,ds.$$

It follows immediately that

(4.12) $$\lim_{t \to 0} \frac{1}{t}[p(r, \psi(t,u)) - p(r,u)] = \frac{1}{r}(\psi(r,u) - u).$$

For $t > 0$ and $u \in U_1$ let

$$q(r,t,u) := [p(r, \psi(t,u)) - p(r,u)](\partial p(r,u))^{-1}.$$

Then we have

(4.13) $$\lim_{t \to 0} \frac{1}{t} q(r,t,u) = \frac{1}{r}(\psi(r,u) - u)(\partial p(r,u))^{-1}.$$

By Taylor's expansion,

$$p(r, \psi(t,u)) - p(r,u) = (\psi(t,u) - u)(\partial p(r,u)) + o(|\psi(t,u) - u|)$$

and consequently

(4.14) $$q(r,t,u) = (\psi(t,u) - u) + o(|\psi(t,u) - u|)$$

as $t \to 0$. It follows that

$$|q(r,t,u)| = |\psi(t,u) - u| + o(|\psi(t,u) - u|)$$

and hence

$$\lim_{t \to 0} |q(r,t,u)|\,|\psi(t,u) - u|^{-1} = 1,$$

which together with (4.13) implies that

$$\lim_{t \to 0} \frac{1}{t}|\psi(t,u) - u| = \lim_{t \to 0} \frac{1}{t}|q(r,t,u)| = \frac{1}{r}|(\psi(r,u) - u)(\partial p(r,u))^{-1}|.$$

Then we get from (4.14) that

(4.15) $$q(r,t,u) = (\psi(t,u) - u) + o(t).$$

From (4.13) and (4.15) it follows that

$$\lim_{t \to 0} \frac{1}{t}(\psi(t,u) - u) = \frac{1}{r}(\psi(r,u) - u)(\partial p(r,u))^{-1}$$

for all $u \in U_1$. Clearly, the entries of $\partial p(r,u)$ and hence those of $(\partial p(r,u))^{-1}$ are continuous in $u \in U_1$. Then the derivative $\psi'_t(0,u)$ exists and is continuous in $u \in U_1$. Since $U_1$ can be arbitrary, we get the desired regularity. □

Now let $(\gamma(t))_{t \geq 0}$ be a stochastically continuous SC-semigroup associated with $(Q(t))_{t \geq 0}$ with characteristic function determined by (3.3) and (3.4). We consider the following.

HYPOTHESIS 4.3. *Suppose that*

$$(4.16) \qquad \int_D [\xi_1 + l_{12}(\xi_2)] m(t, d\xi) < \infty$$

*or, equivalently,*

$$(4.17) \qquad \int_D |\xi| \gamma(t, d\xi) < \infty$$

*for all $t \geq 0$.*

Under the above hypothesis, there is a more convenient representation for the function $\phi(t,u)$. Indeed, from (3.4) it follows that

$$(4.18) \qquad h_1(t) := \frac{\partial \phi}{\partial u_1}(t,u) \bigg|_{u=0} = b_1(t) + \int_D \xi_1 m(t, d\xi)$$

and

$$(4.19) \qquad h_2(t) := \frac{\partial \phi}{\partial u_2}(t,u) \bigg|_{u=0} = b_2(t) + \int_D [\xi_2 - \chi(\xi_2)] m(t, d\xi).$$

By differentiating both sides of (3.3) we get

$$(4.20) \qquad \int_D \xi_1 \gamma(t, d\xi) = h_1(t) \quad \text{and} \quad \int_D \xi_2 \gamma(t, d\xi) = h_2(t).$$

The proof of the next proposition is similar to that of Proposition 4.1.

PROPOSITION 4.2. *If Hypotheses 4.1 and 4.3 hold, the function $\phi(t,u)$ has representation*

$$(4.21) \qquad \begin{aligned} \phi(t,u) &= b_1(t)u_1 + h_2(t)u_2 \\ &\quad + a(t)u_2^2 + \int_D (e^{\langle u, \xi \rangle} - 1 - u_2 \xi_2) m(t, d\xi), \end{aligned}$$

*where $a(t)$, $b_1(t)$ and $m(t, d\xi)$ are as in Proposition 3.1 and $h_2(t)$ is defined by (4.19) and satisfies*

$$(4.22) \quad h_2(r+t) = h_1(r)q_{12}(t) + h_2(r)\beta_{22}(t) + h_2(t), \qquad r, t \geq 0.$$

PROPOSITION 4.3. *If Hypotheses 4.2 and 4.3 hold, the function $t \mapsto h_1(t)$ is continuously differentiable on $[0, \infty)$.*

PROOF. By (3.11), for any $u_1 \in \mathbb{C}_-$ we have

$$\int_D \exp\{u_1 \xi_1\} P(t, x, d\xi) \tag{4.23}$$
$$= \exp\{x_1 \psi_1(t, (u_1, 0)) + \phi(t, (u_1, 0))\}, \qquad u_1 \in \mathbb{C}_-.$$

Then the projection of $P(t, x, \cdot)$ to $\mathbb{R}_+$ is independent of $x_2 \in \mathbb{R}$. Let $P_1(t, x_1, \cdot)$ denote this projection. In view of (4.23), we see that $(P_1(t))_{t \geq 0}$ is the transition semigroup of a CBI-process in the sense of Kawazu and Watanabe [16]. By [16], Theorem 1.1, we have the representation

$$\phi(t, (u_1, 0)) = \int_0^t \bigg[ b_1 \psi_1(s, (u_1, 0)) \tag{4.24}$$
$$+ \int_0^\infty (\exp\{\psi_1(s, (u_1, 0))\xi_1\} - 1) m_1(d\xi_1) \bigg] ds,$$

where $b_1 \geq 0$ is a constant and $m_1(d\xi_1)$ is a $\sigma$-finite measure on $(0, \infty)$ such that

$$\int_0^\infty \chi(\xi_1) m_1(d\xi_1) < \infty.$$

By differentiating both sides of (4.24) with respect to $u_1$ at zero and appealing to (4.3) and (4.18) it is easy to show that

$$h_1(t) = \bigg( b_1 + \int_D \xi_1 m_1(d\xi_1) \bigg) \int_0^t q_{11}(s) \, ds. \tag{4.25}$$

Then we must have

$$\int_0^\infty \xi_1 m_1(d\xi_1) < \infty.$$

In view of (4.5), Hypothesis 4.2 implies that $q_{11}(s)$ is continuous in $s \geq 0$ with $q_{11}(0) = 1$. By (4.25) we find that $h_1(t)$ is differentiable in $t \geq 0$. □

THEOREM 4.2. *Suppose that Hypotheses 4.1 and 4.3 hold. Then $(\gamma(t))_{t \geq 0}$ is regular if and only if $t \mapsto h_2(t)$ is absolutely continuous on $[0, \infty)$.*

PROOF. Based on representation (4.21), this follows as in the proof of Theorem 3.1. □

COROLLARY 4.1. *Suppose that Hypotheses 4.2 and 4.3 hold. Then $(\gamma(t))_{t \geq 0}$ is regular if and only if the function $h_2(\cdot)$ is differentiable at some and hence all $t \geq 0$.*

PROOF. By Proposition 4.3, the function $t \mapsto h_1(t)$ is continuously differentiable. From (4.22) we know that the differentiability of $h_2(\cdot)$ at any $t_0 \geq 0$ implies its differentiability at 0. Using relation (4.22) once again we see that $h_2(\cdot)$ has right derivative at every $t \geq 0$ with

$$h_2'(t+) = h_1'(0+)q_{12}(t) + h_2'(0+)\beta_{22}(t).$$

This function is continuous in $t \geq 0$, so $h_2(\cdot)$ is absolutely continuous. Then the desired result follows from Theorem 4.2. □

**5. One-dimensional stochastic equations.** In this section we prove the existence and pathwise uniqueness of solution of a one-dimensional stochastic equation with non-Lipschitz coefficients and jumps of Poisson type. To simplify the calculations, we only consider a special case for the coefficients which is sufficient for the applications in the next section. The result may be regarded as an extension of the well-known result of Yamada and Watanabe [29]; see also [13] and the references therein for various generalizations of their result in the setting of diffusion processes. For the general background and notation of stochastic equations we refer to [15].

Let $\theta_0 \geq 0$ and $\theta_1 \geq 0$ be constants and let $m(d\xi)$ and $\mu(d\xi)$ be $\sigma$-finite measures on $(0, \infty)$ satisfying $m(l_1) + \mu(l_{12}) < \infty$. Suppose that $(\Omega, \mathscr{F}, \mathscr{F}_t, \mathbf{P})$ is a filtered probability space satisfying the usual hypotheses on which the following are defined:

(a) an $r$-dimensional Brownian motion $B(\cdot) = (B_1(\cdot), \ldots, B_r(\cdot))$;

(b) a Poisson random measure $N_0(ds, d\xi)$ on $(0, \infty)^2$ with intensity $ds\, m(d\xi)$;

(c) a Poisson random measure $N_1(ds, du, d\xi)$ on $(0, \infty)^3$ with intensity $ds\, du\, \mu(d\xi)$;

(d) an $r$-dimensional progressive process $\sigma(\cdot) = (\sigma_1(\cdot), \ldots, \sigma_r(\cdot))$ such that $|\sigma(t)| \leq \bar{\sigma}(t)$ for all $t \geq 0$ and a nonnegative deterministic increasing function $\bar{\sigma}(\cdot)$ on $[0, \infty)$;

(e) a nonnegative progressive process $b(\cdot)$ such that $b(t) \leq \bar{b}(t)$ for all $t \geq 0$ and a nonnegative deterministic increasing function $\bar{b}(\cdot)$ on $[0, \infty)$;

(f) a progressive process $\beta(\cdot)$ such that $|\beta(t)| \leq \bar{\beta}(t)$ for all $t \geq 0$ and a nonnegative deterministic increasing function $\bar{\beta}(\cdot)$ on $[0, \infty)$;

(g) a nonnegative progressive process $l(\cdot)$ such that $l(t) \leq \bar{l}(t)$ for all $t \geq 0$ and a nonnegative deterministic increasing function $\bar{l}(\cdot)$ on $[0, \infty)$.

We assume that $B(\cdot)$, $N_0$ and $N_1$ are independent of each other and $(\mathscr{F}_t)_{t \geq 0}$ is the augmented natural filtration generated by them. Consider the stochastic integral equation

$$x(t) = x(0) + \int_0^t (b(s) + \beta(s)x(s))\, ds$$

$$(5.1) \qquad + \sum_{j=1}^{r} \int_0^t \sigma_j(s)\sqrt{2x(s)}\,dB_j(s) + \int_0^t \int_0^\infty \theta_0 \xi N_0(ds, d\xi)$$

$$+ \int_0^t \int_0^{l(s)x(s-)} \int_0^\infty \theta_1 \xi \tilde{N}_1(ds, du, d\xi),$$

where $\tilde{N}_1(ds, du, d\xi) = N_1(ds, du, d\xi) - ds\,du\,\mu(d\xi)$. Clearly, the diffusion coefficients of (5.1) do not meet the requirements of [15], page 265. Observe also that integration in the last term on the right-hand side is taken over a random set. By a *solution* of (5.1) we mean a nonnegative càdlàg progressive process $x(\cdot)$ satisfying the equation a.s. for each $t \geq 0$. We say *pathwise uniqueness* of solution holds for (5.1) if any two solutions of the equation with the same initial state are indistinguishable.

PROPOSITION 5.1. *Let $x(\cdot)$ be a solution of (5.1) satisfying $\mathbf{E}[x(0)] < \infty$. Then we have*

$$(5.2) \qquad \mathbf{E}[x(t)] \leq \{\mathbf{E}[x(0)] + t\bar{b}(t) + \theta_0 m(l_1) t\} \exp\{t\bar{\beta}(t)\}$$

*for all $t \geq 0$.*

PROOF. Let $\tau_n = \inf\{t \geq 0 : x(t) \geq n\}$ and $x_n(t) = x(t \wedge \tau_n)$. By [12], page 131, we have $\mathbf{E}[x_n(s)] = \mathbf{E}[x_n(s-)] \leq \mathbf{E}[n \vee x(0)] < \infty$ for almost all $s \geq 0$. Since $\bar{b}(t)$ and $\bar{\beta}(t)$ are both increasing in $t \geq 0$, it follows from (5.1) that

$$\mathbf{E}[x_n(t)] \leq \mathbf{E}[x(0)] + t\bar{b}(t) + \theta_0 m(l_1) t + \bar{\beta}(t) \int_0^t \mathbf{E}[x_n(s)]\,ds.$$

By Gronwall's inequality we get

$$\mathbf{E}[x_n(t)] \leq \{\mathbf{E}[x(0)] + t\bar{b}(t) + \theta_0 m(l_1) t\} \exp\{\bar{\beta}(t) t\}.$$

Then (5.2) follows by Fatou's lemma. □

PROPOSITION 5.2. *Let $x(\cdot)$ be a solution of (5.1) satisfying $\mathbf{E}[x(0)] < \infty$. Then we have*

$$\mathbf{E}\left[\sup_{0 \leq s \leq t} x(s)\right] \leq \mathbf{E}[x(0)] + [\bar{b}(t) + \theta_0 m(l_1)] t$$

$$(5.3) \qquad\qquad + [\bar{\beta}(t) + \theta_1 \bar{l}(t)\mu(l_{12})] \int_0^t \mathbf{E}[x(s)]\,ds$$

$$+ 4[2r\bar{\sigma}(t) + \theta_1 \sqrt{\bar{l}(t)}\mu(l_{12})]\left(1 + \int_0^t \mathbf{E}[x(s)]\,ds\right)$$

*for all $t \geq 0$.*

PROOF. Applying Doob's inequality to the martingale terms in (5.1),

$$\mathbf{E}\bigg[\sup_{0\le s\le t} x(s)\bigg] \le \mathbf{E}[x(0)] + \int_0^t \mathbf{E}[\bar{b}(s) + \bar{\beta}(s)x(s)]\,ds + \theta_0 t \int_0^\infty \xi m(d\xi)$$

$$+ 4\sum_{j=1}^r \mathbf{E}^{1/2}\bigg[\bigg(\int_0^t \sigma_j(s)\sqrt{2x(s)}\,dB_j(s)\bigg)^2\bigg]$$

$$+ 4\mathbf{E}^{1/2}\bigg[\bigg(\int_0^t \int_0^{l(s)x(s-)} \int_0^1 \theta_1 \xi \tilde{N}_1(ds,du,d\xi)\bigg)^2\bigg]$$

$$+ \mathbf{E}\bigg[\int_0^t \int_0^{l(s)x(s-)} \int_1^\infty \theta_1 \xi N_1(ds,du,d\xi)\bigg]$$

$$\le \mathbf{E}[x(0)] + \int_0^t \bar{b}(s)\,ds + \int_0^t \bar{\beta}(s)\mathbf{E}[x(s)]\,ds + \theta_0 t m(l_1)$$

$$+ 8\sum_{j=1}^r \bigg(\int_0^t \mathbf{E}[\sigma_j^2(s)x(s)]\,ds\bigg)^{1/2}$$

$$+ 4\theta_1 \mu(l_{12})\bigg(\int_0^t \mathbf{E}[l(s)x(s)]\,ds\bigg)^{1/2}$$

$$+ \theta_1 \mu(l_{12})\int_0^t \mathbf{E}[l(s)x(s)]\,ds.$$

Then we obtain (5.3) by combining the terms. $\square$

THEOREM 5.1. *Suppose that $x_1(\cdot)$ and $x_2(\cdot)$ are two solutions of* (5.1) *satisfying* $\mathbf{E}[x_2(0) + x_1(0)] < \infty$. *Then we have*

(5.4) $$\mathbf{E}[|x_2(t) - x_1(t)|] \le \mathbf{E}[|x_2(0) - x_1(0)|]\exp\{t\bar{\beta}(t)\}.$$

*Consequently, the pathwise uniqueness of solution holds for* (5.1).

PROOF. By Proposition 5.1 it is easy to find that $\mathbf{E}[x_2(t) + x_1(t)]$ is locally bounded in $t \ge 0$. Let $z(t) = x_2(t) - x_1(t)$. Then we have

(5.5)
$$z(t) = z(0) + \int_0^t \beta(s)z(s)\,ds$$
$$+ \sum_{j=1}^r \int_0^t \sigma_j(s)(\sqrt{2x_2(s)} - \sqrt{2x_1(s)})\,dB_j(s)$$
$$+ \int_0^t \int_{l(s)x_1(s-)}^{l(s)x_2(s-)} \int_0^\infty \theta_1 \xi \mathbb{1}_{\{x_1(s-)\le x_2(s-)\}} \tilde{N}_1(ds,du,d\xi)$$
$$- \int_0^t \int_{l(s)x_2(s-)}^{l(s)x_1(s-)} \int_0^\infty \theta_1 \xi \mathbb{1}_{\{x_1(s-)>x_2(s-)\}} \tilde{N}_1(ds,du,d\xi).$$

Let $\{a_k\}$ be the sequence defined inductively by $a_0 = 1$ and $a_k = a_{k-1}e^{-k}$ for $k \geq 1$. It is easy to check that $\int_{a_k}^{a_{k-1}} (ku)^{-1} du = 1$. For $k \geq 1$ let $g_k$ be a non-negative continuous function on $[0, \infty)$ which has support contained in $[a_k, a_{k-1}]$ and satisfies $0 \leq g_k(u) \leq 2(ku)^{-1}$ and $\int_{a_k}^{a_{k-1}} g_k(u) \, du = 1$. Then

$$h_k(x) = \int_0^{|x|} dy \int_0^y g_k(u) \, du, \qquad x \in \mathbb{R},$$

defines a twice continuously differentiable function $h_k$ such that $h_k(x) \to |x|$ increasingly as $k \to \infty$. Set $H_k(x, \xi) = \Delta_\xi h_k(x) - h_k'(x)\xi$. By (5.5) and Itô's formula,

(5.6)
$$\begin{aligned}
h_k(z(t)) &= h_k(z(0)) + \int_0^t h_k'(z(s))\beta(s)z(s) \, ds + \text{ mart.} \\
&+ \sum_{j=1}^r \int_0^t h_k''(z(s))\sigma_j^2(s)(\sqrt{x_2(s)} - \sqrt{x_1(s)})^2 \, ds \\
&+ \int_0^t ds \int_0^\infty H_k(z(s), \theta_1\xi)l(s)z(s)\mathbb{1}_{\{z(s) \geq 0\}}\mu(d\xi) \\
&- \int_0^t ds \int_0^\infty H_k(z(s), -\theta_1\xi)l(s)z(s)\mathbb{1}_{\{z(s) < 0\}}\mu(d\xi);
\end{aligned}$$

see, for example, [9], pages 334–335. Note that $|h_k'(x)| \leq 1$ and $0 \leq h_k''(x) = g_k(|x|) \leq 2|kx|^{-1}$. It follows that

$$h_k''(z(s))(\sqrt{x_2(s)} - \sqrt{x_1(s)})^2 \leq h_k''(z(s))|x_2(s) - x_1(s)| \leq 2/k.$$

By the mean-value theorem and Taylor's expansion it is easy to show that $|H_k(x, \xi)x| \leq |2\xi x| \wedge |k^{-1}\xi^2|$ whenever $x\xi \geq 0$. Then we may take the expectations in (5.6) to find

(5.7)
$$\begin{aligned}
\mathbf{E}[h_k(z(t))] &\leq \mathbf{E}[h_k(z(0))] + \bar{\beta}(t) \int_0^t \mathbf{E}[|z(s)|] \, ds \\
&+ \frac{2}{k} \sum_{j=1}^r \int_0^t \mathbf{E}[\sigma_j^2(s)] \, ds \\
&+ \theta_1 \bar{l}(t) \int_0^t ds \int_0^\infty \mathbf{E}[2\xi|z(s)|] \wedge (k^{-1}\theta_1\xi^2)\mu(d\xi).
\end{aligned}$$

Letting $k \to \infty$ in (5.7) we obtain

$$\mathbf{E}[|z(t)|] \leq \mathbf{E}[|z(0)|] + \bar{\beta}(t) \int_0^t \mathbf{E}[|z(s)|] \, ds.$$

Then (5.4) follows by Gronwall's inequality. $\square$

Now we turn to the existence of the solution of (5.1). The Picard iteration method fails for this equation because the diffusion coefficients are not Lipschitz. Since the coefficients are random, we cannot follow the standard argument of martingale problem. In the approach given below, we first approximate the random coefficients by some simple processes and consider a sequence of equations without small and large jumps. The original coefficients and the small jumps are retrieved by a limit argument based on the second moment analysis. Finally, we obtain the solution of (5.1) by adding the large jumps.

A stochastic process $q(\cdot)$ defined on $(\Omega, \mathscr{F}, \mathscr{F}_t, \mathbf{P})$ is called a *simple process* if there is a sequence $0 = r_0 < r_1 < r_2 < \cdots$ increasing to infinity and a sequence of random variables $\{\eta_k\}$ such that $\eta_k$ is $\mathscr{F}_{r_k}$-measurable and

$$(5.8) \qquad q(t) = \eta_0 \mathbb{1}_{\{0\}}(t) + \sum_{k=0}^{\infty} \eta_k \mathbb{1}_{(r_k, r_{k+1}]}(t), \qquad t \geq 0.$$

We approximate the coefficients of (5.1) in the following way:

(a) Let $\{\sigma_n\}$ be a sequence of $r$-dimensional simple processes such that $|\sigma_n(t)| \leq \bar{\sigma}(t)$ for all $t \geq 0$ and $\sigma_n(\cdot) \to \sigma(\cdot)$ a.s. in $L^2([0, J], \lambda)$ for all integers $J \geq 1$.

(b) Let $\{b_n\}$ be a sequence of nonnegative simple processes such that $b_n(t) \leq \bar{b}(t)$ for all $t \geq 0$ and $b_n(\cdot) \to b(\cdot)$ a.s. in $L^2([0, J], \lambda)$ for all integers $J \geq 1$.

(c) Let $\{\beta_n\}$ be a sequence of simple processes such that $|\beta_n(t)| \leq \bar{\beta}(t)$ for all $t \geq 0$ and $\beta_n(\cdot) \to \beta(\cdot)$ a.s. in $L^2([0, J], \lambda)$ for all integers $J \geq 1$.

(d) Let $\{l_n\}$ be a sequence of simple processes such that $l_n(t) \leq \bar{l}(t)$ for all $t \geq 0$ and $l_n(\cdot) \to l(\cdot)$ a.s. in $L^2([0, J], \lambda)$ for all integers $J \geq 1$.

Let $L \geq 1$ be an integer and let $\{\varepsilon_n\}$ be a decreasing sequence such that $\mu(\{\varepsilon_n : n \geq 1\}) = 0$, $0 < \varepsilon_n \leq 1$ and $\varepsilon_n \to 0$ as $n \to \infty$. Suppose that $x(0)$ is a nonnegative $\mathscr{F}_0$-measurable random variable satisfying $\mathbf{E}[x(0)] < \infty$. Let $x_n(\cdot)$ denote the nonnegative solution of the stochastic equation

$$(5.9) \qquad \begin{aligned} x_n(t) &= x(0) + \int_0^t (b_n(s) + \beta_n(s) x_n(s))\, ds \\ &\quad + \sum_{j=1}^r \int_0^t \sigma_{n,j}(s)\sqrt{2 x_n(s)}\, dB_j(s) + \int_0^t \int_{\varepsilon_n}^L \theta_0 \xi N_0(ds, d\xi) \\ &\quad + \int_0^t \int_0^{l_n(s) x_n(s-)} \int_{\varepsilon_n}^L \theta_1 \xi \tilde{N}_1(ds, du, d\xi). \end{aligned}$$

Based on Proposition 5.2 and the results in [15], pages 235–237, the existence of the strong solution of the above equation follows by arguments similar to those of [15], pages 245–246. Let

$$y_{n,j}(t) := \int_0^t \sigma_{n,j}(s)\sqrt{2x_n(s)}\,dB_j(s) \tag{5.10}$$

and

$$z_n(t) := \int_0^t \int_0^{l_n(s)x_n(s-)} \int_{\varepsilon_n}^L \xi \tilde{N}_1(ds,du,d\xi). \tag{5.11}$$

LEMMA 5.1. *For $1 \leq j \leq r$ the sequence $y_{n,j}(\cdot)$ is tight in $C([0,\infty),\mathbb{R})$, and the sequences $x_n(\cdot)$ and $z_n(\cdot)$ are tight in $D([0,\infty),\mathbb{R})$.*

PROOF. By Proposition 5.1 it is easy to show that $C(t) := \sup_{n \geq 1} \mathbf{E}[x_n(t)]$ is a locally bounded function of $t \geq 0$. By (5.10) we have

$$\mathbf{E}[|y_{n,j}(t)|^2] = 2\int_0^t \mathbf{E}[\sigma_{n,j}^2(s)x_n(s)]\,ds \leq 2\int_0^t \bar{\sigma}^2(s)C(s)\,ds \tag{5.12}$$

and

$$\mathbf{E}[|z_n(t)|^2] = \int_0^t \mathbf{E}[l_n(s)x_n(s)]\,ds \int_{\varepsilon_n}^L \xi^2\mu(d\xi)$$
$$\leq \int_0^t \bar{l}(s)C(s)\,ds \int_0^L \xi^2\mu(d\xi). \tag{5.13}$$

Then $y_{n,j}(t)$ and $z_n(t)$ are tight sequences of random variables for every fixed $t \geq 0$. Now let $\{\tau_n\}$ be a sequence of stopping times bounded above by $T \geq 0$. By the properties of independent increments of the Brownian motion and the Poisson process we obtain as in the calculations in (5.12) and (5.13) that

$$\mathbf{E}[|y_{n,j}(\tau_n + t) - y_{n,j}(\tau_n)|^2] \leq 2\int_0^t \bar{\sigma}^2(T+s)C(T+s)\,ds$$

and

$$\mathbf{E}[|z_n(\tau_n + t) - z_n(\tau_n)|^2] \leq \int_0^t \bar{l}(T+s)C(T+s)\,ds \int_0^L \xi^2\mu(d\xi).$$

Then $y_{n,j}(\cdot)$ and $z_n(\cdot)$ are tight in $D([0,\infty),\mathbb{R})$ by the criterion of Aldous [1]. Since $C([0,\infty),\mathbb{R})$ is a closed subset of $D([0,\infty),\mathbb{R})$, we infer that $y_{n,j}(\cdot)$ is also tight in $C([0,\infty),\mathbb{R})$. By similar calculations for other terms on the right-hand side of (5.9) we find that $x_n(\cdot)$ is tight in $D([0,\infty),\mathbb{R}_+)$. □

By Lemma 5.1 we may construct a new filtered probability space $(\Omega, \mathscr{F}, \mathscr{F}_t, \mathbf{P})$ satisfying the usual hypotheses on which the following stochastic

equations are realized:

$$x_n(t) = x_n(0) + \int_0^t (b_n(s) + \beta_n(s)x_n(s))\,ds$$

(5.14)
$$+ \sum_{j=1}^r \int_0^t \sigma_{n,j}(s)\sqrt{2x_n(s)}\,dB_{n,j}(s) + \int_0^t \int_{\varepsilon_n}^L \theta_0 \xi N_{n,0}(ds, d\xi)$$

$$+ \int_0^t \int_0^{l_n(s)x_n(s-)} \int_{\varepsilon_n}^L \theta_1 \xi \tilde{N}_{n,1}(ds, du, d\xi),$$

where the processes $\{x_n, B_n, \sigma_n, b_n, \beta_n, l_n\}$ and the random measures $\{N_{n,0}, N_{n,1}\}$ are distributed as $\{x_n, B, \sigma_n, b_n, \beta_n, l_n\}$ and $\{N_0, N_1\}$ in (5.9). Moreover, as $n \to \infty$ we have:

(a) $x_n(\cdot) \to$ a process $x(\cdot)$ a.s. by the topology of $D([0,\infty), \mathbb{R}_+)$;

(b) $B_n(\cdot) \to$ an $r$-dimensional Brownian motion $B(\cdot)$ a.s. by the topology of $C([0,\infty), \mathbb{R}^r)$;

(c) $\xi N_{n,0}(ds, d\xi) \to \xi N_0(ds, d\xi)$ a.s. by the weak convergence of finite measures on $(0, J] \times (0, L]$ for all integers $J \geq 1$, where $N_0(ds, dy)$ is a Poisson random measure on $(0, \infty)^2$ with intensity $ds\,m(dy)$;

(d) $\xi^2 N_{n,1}(ds, du, d\xi) \to \xi^2 N_1(ds, du, d\xi)$ a.s. by the week convergence of finite measures on $(0, J]^2 \times (0, L]$ for all integers $J \geq 1$, where $N_1(ds, du, dy)$ is a Poisson random measure on $(0, \infty)^3$ with intensity $ds\,du\,\mu(dy)$;

(e) $\sigma_n(\cdot)$, $b_n(\cdot)$, $\beta_n(\cdot)$ and $l_n(\cdot)$ converge a.s. to processes $\sigma(\cdot)$, $b(\cdot)$, $\beta(\cdot)$ and $l(\cdot)$, respectively, by the topology of $L^2([0, J], \lambda)$ for each integer $J \geq 1$;

(f) for each $1 \leq j \leq r$ it holds that

(5.15) $$y_{n,j}(t) := \int_0^t \sigma_{n,j}(s)\sqrt{2x_n(s)}\,dB_{n,j}(s) \to \text{ a process } y_j(t) \quad \text{a.s.}$$

by the topology of $C([0,\infty), \mathbb{R})$;

(g) it holds that

(5.16)
$$z_n(t) := \int_0^t \int_0^{l_n(s)x_n(s-)} \int_{\varepsilon_n}^L \xi \tilde{N}_{n,1}(ds, du, d\xi)$$
$$\to \text{ a process } z(t) \quad \text{a.s.}$$

by the topology of $D([0,\infty), \mathbb{R})$.

The existence of such a probability space follows by the Skorokhod representation; see, for example, [12], page 102. Indeed, we can and do assume that the probability space is constructed so that the above assertions hold simultaneously for all integers $L \geq 1$. Of course, the processes $\{x_n(\cdot), y_{n,j}(\cdot), z_n(\cdot), x(\cdot), y_j(\cdot), z(\cdot)\}$ all depend on $L \geq 1$. We suppress this dependence for simplicity of the notation. Note also that the processes $\{B(\cdot), \sigma(\cdot), b(\cdot), \beta(\cdot), l(\cdot)\}$ and the random measures $\{N_0(ds, d\xi), N_1(ds, du, d\xi)\}$ are distributed as those in (5.1).

LEMMA 5.2. *For each $t \geq 0$ we have a.s.*

$$y_j(t) = \int_0^t \sigma_j(s)\sqrt{2x(s)}\, dB_j(s) \tag{5.17}$$

*and*

$$z(t) = \int_0^t \int_0^{l(s)x(s-)} \int_0^L \xi \tilde{N}_{1,n}(ds, du, d\xi). \tag{5.18}$$

PROOF. For $m \geq 1$ let $\tau_m = \inf\{t \geq 0 : x(t) \geq m \text{ or } x_n(t) \geq m \text{ for some } n \geq 1\}$. Let $\{q_{k,m}(\cdot)\}$ be a sequence of nonnegative simple processes bounded above by $m$ such that $q_{k,m}(s) \to x(s)\mathbb{1}_{\{s \leq \tau_m\}}$ a.s. by the topology of $L^2([0,T], \lambda)$ for each $T \geq 0$. By (5.15) we have

$$y_{n,j}(t \wedge \tau_m) = \int_0^t \sigma_{k,j}(s)\sqrt{2q_{k,m}(s)}\, dB_{n,j}(s) + \eta_{n,k,m,j}(t), \tag{5.19}$$

where

$$\eta_{n,k,m,j}(t) = \int_0^t \sigma_{n,j}(s)(\sqrt{2x_n(s)} - \sqrt{2x(s)})\mathbb{1}_{\{s \leq \tau_m\}}\, dB_{n,j}(s)$$
$$+ \int_0^t \sigma_{n,j}(s)(\sqrt{2x(s)} - \sqrt{2q_{k,m}(s)})\mathbb{1}_{\{s \leq \tau_m\}}\, dB_{n,j}(s)$$
$$+ \int_0^t (\sigma_{n,j}(s) - \sigma_{k,j}(s))\sqrt{2q_{k,m}(s)}\, dB_{n,j}(s).$$

It is simple to see that

$$\mathbf{E}[\eta_{n,k,m,j}^2(t)] \leq 6\bar{\sigma}(t)^2 \int_0^t \mathbf{E}[(\sqrt{x_n(s)} - \sqrt{x(s)})^2 \mathbb{1}_{\{s \leq \tau_m\}}]\, ds$$
$$+ 6\bar{\sigma}(t)^2 \int_0^t \mathbf{E}[(\sqrt{x(s)} - \sqrt{q_{k,m}(s)})^2 \mathbb{1}_{\{s \leq \tau_m\}}]\, ds \tag{5.20}$$
$$+ 6m \int_0^t \mathbf{E}[(\sigma_{n,j}(s) - \sigma_{k,j}(s))^2]\, ds.$$

In view of (5.19), the limit $\eta_{k,m,j}(t) = \lim_{n \to \infty} \eta_{n,k,m,j}(t)$ exists and

$$y_j(t \wedge \tau_m) = \int_0^t \sigma_{k,j}(s)\sqrt{2q_{k,m}(s)}\, dB_j(s) + \eta_{k,m,j}(t). \tag{5.21}$$

By (5.20) and Fatou's lemma,

$$\mathbf{E}[\eta_{k,m,j}^2(t)] \leq 6\bar{\sigma}(t)^2 \int_0^t \mathbf{E}[(\sqrt{x(s)} - \sqrt{q_{k,m}(s)})^2 \mathbb{1}_{\{s \leq \tau_m\}}]\, ds$$
$$+ 6m \int_0^t \mathbf{E}[(\sigma_j(s) - \sigma_{k,j}(s))^2]\, ds,$$

which goes to zero as $k \to \infty$. Now we can take the limit in (5.21) to obtain

$$y_j(t \wedge \tau_m) = \int_0^t \sigma_j(s)\sqrt{2x(s)}\mathbb{1}_{\{s \le \tau_m\}}\, dB_j(s).$$

Then we have (5.17) since $\tau_m \to \infty$ as $m \to \infty$. Equality (5.18) can be proved using similar ideas. □

LEMMA 5.3. *For each $t \ge 0$ we have a.s.*

(5.22)
$$\begin{aligned}
x(t) = {} & x(0) + \int_0^t (b(s) + \beta(s)x(s))\, ds \\
& + \sum_{j=1}^r \int_0^t \sigma_j(s)\sqrt{2x(s)}\, dB_j(s) + \int_0^t \int_0^L \theta_0 \xi N_0(ds, d\xi) \\
& + \int_0^t \int_0^{l(s)x(s-)} \int_0^L \theta_1 \xi \tilde{N}_1(ds, du, d\xi).
\end{aligned}$$

PROOF. By dominated convergence we have a.s.

$$\lim_{n\to\infty} \int_0^t (b_n(s) + \beta_n(s)x_n(s))\, ds = \int_0^t (b(s) + \beta(s)x(s))\, ds.$$

On the other hand, it is easy to show that a.s.

$$\lim_{n\to\infty} \int_0^t \int_{\varepsilon_n}^L \xi N_{n,0}(ds, d\xi) = \int_0^t \int_0^L \xi N_0(ds, d\xi).$$

Then (5.22) follows from (5.14) and Lemma 5.2. □

THEOREM 5.2. *There is a solution $x(\cdot)$ of (5.1).*

PROOF. By Lemma 5.3 there is a sequence of processes $\{x_k(\cdot)\}$ satisfying the equations

$$\begin{aligned}
x_k(t) = {} & x(0) + \int_0^t (b(s) + \beta(s)x_k(s))\, ds + \sum_{j=1}^r \int_0^t \sigma_j(s)\sqrt{2x_k(s)}\, dB_j(s) \\
& + \int_0^t \int_0^k \theta_0 \xi N_0(ds, d\xi) + \int_0^t \int_0^{l(s)x_k(s-)} \int_0^1 \theta_1 \xi \tilde{N}_1(ds, du, d\xi) \\
& + \int_0^t \int_0^{l(s)x_k(s-)} \int_1^k \theta_1 \xi N_1(ds, du, d\xi) - \int_1^\infty \theta_1 \xi \mu(d\xi) \int_0^t l(s)x_k(s)\, ds.
\end{aligned}$$

The pathwise uniqueness of solutions holds for those equations by Theorem 5.1. Based on this fact, it is easy to show that $x_k(\cdot)$ is increasing in $k \ge 1$. Let $x(\cdot) := \lim_{k\to\infty} x_k(\cdot)$. By Propositions 5.1 and 5.2 and Fatou's lemma we conclude that $\mathbf{E}[\sup_{0\leq s\leq T} x(s)] < \infty$ for each $T \geq 0$. Then we infer that $x(\cdot)$ satisfies (5.1). □

In particular, if $\{\sigma, b, \beta, l\}$ are all deterministic constants, Theorems 5.1 and 5.2 imply that (5.1) has a unique *strong solution* $x(\cdot)$ and the solution is a strong Markov process; see, for example, [15], pages 163–166 and page 215. By Itô's formula, we find that $x(\cdot)$ has generator $L$ determined by

$$Lf(x) = \alpha x f''(x) + (b + \beta x)f'(x) + \int_0^\infty \Delta_{\theta_0 \xi} f(x) m(d\xi)$$
(5.23)
$$+ \int_0^\infty (\Delta_{\theta_1 \xi} f(x) - f'(x)\theta_1 \xi) l x \mu(d\xi),$$

where $\alpha = \sum_{j=1}^r \sigma_j^2$. Then $x(\cdot)$ is a CBI-process; see [16] and [27]. The stochastic equation (5.1) gives explicit representations of the two types of jumps of the process in terms of the Poisson random measures $N_0(ds, d\xi)$ and $N_1(ds, du, d\xi)$. As far as we know, this characterization of the CBI-process has not appeared in the literature. In the general case, the solution of (5.1) can be regarded as a *generalized CBI-process* with random parameters.

**6. Constructions of the two-dimensional processes.** Based on the results in the last section, we here construct two classes of Markov processes as strong solutions of stochastic integral equations. The first class is the regular affine process and the second is the catalytic CBI-process. The characterizations of those processes in terms of stochastic equations play the key role in the study of the limit theorems in the next section. To simplify the discussions, we impose some conditions on the jumps so that the processes possess finite first moments.

DEFINITION 6.1. A set of parameters $(a, (\alpha_{ij}), (b_1, b_2), (\beta_{ij}), m, \mu)$ is called *admissible* if:

(i) $a \in \mathbb{R}_+$ is a constant;
(ii) $(\alpha_{ij})$ is a symmetric nonnegative definite $(2 \times 2)$-matrix;
(iii) $(b_1, b_2) \in D$ is a vector;
(iv) $(\beta_{ij})$ is a $(2 \times 2)$-matrix with $\beta_{12} = 0$;
(v) $m(d\xi)$ is a $\sigma$-finite measure on $D$ supported by $D \setminus \{0\}$ such that

$$\int_D [l_1(\xi_1) + l_{12}(\xi_2)] m(d\xi) < \infty;$$

(vi) $\mu(d\xi)$ is a $\sigma$-finite measure on $D$ supported by $D \setminus \{0\}$ such that

$$\int_D [l_{12}(\xi_1) + l_{12}(\xi_2)] \mu(d\xi) < \infty.$$

THEOREM 6.1 ([10]). *Suppose that $(a, (\alpha_{ij}), (b_1, b_2), (\beta_{ij}), m, \mu)$ is a set of admissible parameters. For $u = (u_1, u_2) \in U$ set*

$$F(u) = b_1 u_1 + b_2 u_2 + a u_2^2 + \int_D (e^{\langle u, \xi \rangle} - 1 - \xi_2 u_2) m(d\xi) \tag{6.1}$$

*and*

$$\begin{aligned} R(u) &= \beta_{11} u_1 + \beta_{21} u_2 + \alpha_{11} u_1^2 + 2\alpha_{12} u_1 u_2 + \alpha_{22} u_2^2 \\ &\quad + \int_D (e^{\langle u, \xi \rangle} - 1 - \xi_1 u_1 - \xi_2 u_2) \mu(d\xi). \end{aligned} \tag{6.2}$$

*Then there is a unique regular affine semigroup $(P(t))_{t \geq 0}$ determined by (3.11) where $\psi_2(t, u) = e^{\beta_{22} t} u_2$, $\psi_1(t, u)$ solves the generalized Riccati equation*

$$\frac{\partial}{\partial t} \psi_1(t, u) = R(\psi_1(t, u), e^{\beta_{22} t} u_2), \qquad \psi_1(0, u) = u_1 \tag{6.3}$$

*and*

$$\phi(t, u) = \int_0^t F(\psi_1(s, u), e^{\beta_{22} s} u_2) \, ds. \tag{6.4}$$

Let $(a, (\alpha_{ij}), (\beta_{ij}), (b_j), m, \mu)$ be a set of admissible parameters and let $A$ be the generator of the regular affine semigroup $(P(t))_{t \geq 0}$ characterized by Theorem 6.1. It is not hard to show that

$$\begin{aligned} Af(x) &= \alpha_{11} x_1 f_{11}''(x) + 2\alpha_{12} x_1 f_{12}''(x) + \alpha_{22} x_1 f_{22}''(x) + a f_{22}''(x) \\ &\quad + (b_1 + \beta_{11} x_1) f_1'(x) + (b_2 + \beta_{21} x_1 + \beta_{22} x_2) f_2'(x) \\ &\quad + \int_D (\Delta_\xi f(x) - f_2'(x) \xi_2) m(d\xi) \\ &\quad + \int_D (\Delta_\xi f(x) - \langle \nabla f(x), \xi \rangle) x_1 \mu(d\xi) \end{aligned} \tag{6.5}$$

for $f \in C^2(D)$, where $\nabla f(x) = (f_1'(x), f_2'(x))$.

Let $\sigma_0 = \sqrt{a}$ and let $(\sigma_{ij})$ be a $(2 \times 2)$-matrix satisfying $(\alpha_{ij}) = (\sigma_{ij})(\sigma_{ij})^\tau$. Let $(\Omega, \mathscr{F}, \mathscr{F}_t, \mathbf{P})$ be a filtered probability space satisfying the usual hypotheses. Suppose that on this probability space the following objects are defined:

(a) a three-dimensional Brownian motion $B(\cdot) = (B_0(\cdot), B_1(\cdot), B_2(\cdot))$;

(b) a Poisson random measure $N_0(ds, d\xi)$ on $(0, \infty) \times D$ with intensity $ds\, m(d\xi)$;

(c) a Poisson random measure $N_1(ds, du, d\xi)$ on $(0, \infty)^2 \times D$ with intensity $ds\, du\, \mu(d\xi)$.

We assume that $B_0(\cdot)$, $N_0$ and $N_1$ are independent of each other and $(\mathscr{F}_t)_{t\geq 0}$ is the augmented natural filtration generated by them. Let $x(0)$ be a nonnegative $\mathscr{F}_0$-measurable random variable defined on $(\Omega, \mathscr{F}, \mathscr{F}_t, \mathbf{P})$. By Theorems 5.1 and 5.2 there is a unique strong solution $x(\cdot)$ of

$$
\begin{aligned}
x(t) = x(0) &+ \int_0^t (b_1 + \beta_{11} x(s))\, ds + \int_0^t \sigma_{11}\sqrt{2x(s)}\, dB_1(s) \\
&+ \int_0^t \sigma_{12}\sqrt{2x(s)}\, dB_2(s) + \int_0^t \int_D \xi_1 N_0(ds, d\xi) \\
&+ \int_0^t \int_0^{x(s-)} \int_D \xi_1 \tilde{N}_1(ds, du, d\xi).
\end{aligned}
\tag{6.6}
$$

As explained at the end of the last section, $x(\cdot)$ is a CBI-process. In addition, let $z(0)$ be an $\mathscr{F}_0$-measurable random variable defined on $(\Omega, \mathscr{F}, \mathscr{F}_t, \mathbf{P})$. We consider the equation

$$
\begin{aligned}
z(t) = z(0) &+ \int_0^t (b_2 + \beta_{21}x(s) + \beta_{22}z(s))\, ds + \int_0^t \sqrt{2}\sigma_0\, dB_0(s) \\
&+ \int_0^t \sigma_{21}\sqrt{2x(s)}\, dB_1(s) + \int_0^t \sigma_{22}\sqrt{2x(s)}\, dB_2(s) \\
&+ \int_0^t \int_D \xi_2 \tilde{N}_0(ds, d\xi) + \int_0^t \int_0^{x(s-)} \int_D \xi_2 \tilde{N}_1(ds, du, d\xi).
\end{aligned}
\tag{6.7}
$$

THEOREM 6.2.  *The equation system* (6.6) *and* (6.7) *has a unique strong solution* $(x(\cdot), z(\cdot))$. *Moreover*, $(x(\cdot), z(\cdot))$ *is an affine Markov process with generator $A$ given by* (6.5).

PROOF. By Itô's formula it is not hard to show that

$$
\begin{aligned}
z(t) = e^{\beta_{22}t}z(0) &+ e^{\beta_{22}t}\int_0^t e^{-\beta_{22}s}(b_2 + \beta_{21}x(s))\, ds \\
&+ e^{\beta_{22}t}\int_0^t \sqrt{2}\sigma_0 e^{-\beta_{22}s}\, dB_0(s) + e^{\beta_{22}t}\int_0^t \sigma_{21}e^{-\beta_{22}s}\sqrt{2x(s)}\, dB_1(s) \\
&+ e^{\beta_{22}t}\int_0^t \sigma_{22}e^{-\beta_{22}s}\sqrt{2x(s)}\, dB_2(s) + e^{\beta_{22}t}\int_0^t \int_D e^{-\beta_{22}s}\xi_2\tilde{N}_0(ds, d\xi) \\
&+ e^{\beta_{22}t}\int_0^t \int_0^{x(s-)} \int_D e^{-\beta_{22}s}\xi_2\tilde{N}_1(ds, du, d\xi)
\end{aligned}
\tag{6.8}
$$

defines a solution of (6.7) and conversely any solution of (6.7) must be given by (6.8). The uniqueness implies the strong Markov property of $(x(\cdot), z(\cdot))$. By Itô's formula, we find that the Markov process $(x(\cdot), z(\cdot))$ has generator $A$. □

Now suppose that $b_2 \geq 0$ and $m_2(l_1) \leq \infty$, were $m_2$ denotes the projection of $m$ to $\mathbb{R}$. Let $D_+ = \mathbb{R}_+^2 \subset D$. Given a nonnegative $\mathscr{F}_0$-measurable random variable $y(0)$ defined on $(\Omega, \mathscr{F}, \mathscr{F}_t, \mathbf{P})$, we consider the equation

$$\begin{aligned}
y(t) = y(0) &+ \int_0^t (b_2 + \beta_{21} x(s) y(s) + \beta_{22} y(s))\, ds + \int_0^t \sigma_0 \sqrt{2y(s)}\, dB_0(s) \\
&+ \int_0^t \sigma_{21} \sqrt{2x(s)y(s)}\, dB_1(s) + \int_0^t \sigma_{22} \sqrt{2x(s)y(s)}\, dB_2(s) \\
&+ \int_0^t \int_{D_+} \xi_2 N_0(ds, d\xi) + \int_0^t \int_0^{lx(s-)y(s-)} \int_{D_+} \xi_2 \tilde{N}_1(ds, du, d\xi).
\end{aligned} \tag{6.9}$$

A solution $y(\cdot)$ of (6.9) can be regarded as a generalized CBI-process with random parameters governed by the process $x(\cdot)$. Following Dawson and Fleischmann [5], we shall call the pair $(x(\cdot), y(\cdot))$ a *catalytic CBI-process*, where $x(\cdot)$ is the *catalyst process* and $y(\cdot)$ is the *reactant process*.

THEOREM 6.3. *The equation system* (6.6) *and* (6.9) *has a unique strong solution* $(x(\cdot), y(\cdot))$.

PROOF. It suffices to consider the case where the initial states $x(0)$ and $y(0)$ are deterministic. For $n \geq x(0)$ let $\tau_n = \inf\{s \geq 0 : x(s) \geq n\}$ and $x_n(t) = x(t \wedge \tau_n)$. By Theorems 5.1 and 5.2, there is a unique strong solution $(x(\cdot), y_n(\cdot))$ of the equation system formed by (6.6) and

$$\begin{aligned}
y_n(t) = y(0) &+ \int_0^t (b_2 + \beta_{21} x_n(s) y_n(s) + \beta_{22} y_n(s))\, ds \\
&+ \int_0^t \sigma_0 \sqrt{2y_n(s)}\, dB_0(s) + \int_0^t \sigma_{21} \sqrt{2x_n(s)y_n(s)}\, dB_1(s) \\
&+ \int_0^t \sigma_{22} \sqrt{2x_n(s)y_n(s)}\, dB_2(s) + \int_0^t \int_{D_+} \xi_2 N_0(ds, d\xi) \\
&+ \int_0^t \int_0^{lx_n(s-)y_n(s-)} \int_{D_+} \xi_2 \tilde{N}_1(ds, du, d\xi).
\end{aligned} \tag{6.10}$$

By the uniqueness, for any $n \geq m \geq x(0)$ the two processes $y_n(t \wedge \tau_m)$ and $y_m(t \wedge \tau_m)$ are indistinguishable. Since $\tau_n \to \infty$ as $n \to \infty$, it is easy to see that $y(t) := \lim_{n \to \infty} y_n(t)$ is the unique solution of (6.9). $\square$

By Theorem 6.3, the catalytic CBI-process $(x(\cdot), y(\cdot))$ is a strong Markov process with state space $D_+$. Let $D_- = \mathbb{R}_+ \times \mathbb{R}_-$. By Itô's formula we find that $(x(\cdot), y(\cdot))$ has generator $L$ determined by

$$Lf(x) = \alpha_{11} x_1 f''_{11}(x) + 2\alpha_{12} x_1 \sqrt{x_2} f''_{12}(x) + \alpha_{22} x_1 x_2 f''_{22}(x) + a x_2 f''_{22}(x)$$

$$+ (b_1 + \beta_{11}x_1)f_1'(x) + (b_2 + \beta_{21}x_1x_2 + \beta_{22}x_2)f_2'(x)$$

$$+ \int_{D_+} \Delta_\xi f(x) m(d\xi) + \int_{D_-} \Delta_{(\xi_1,0)} f(x) m(d\xi)$$

(6.11)
$$+ \int_{D_+} [\Delta_\xi f(x) - \langle \nabla f(x), \xi \rangle](x_1 \wedge lx_1x_2)\mu(d\xi)$$

$$+ \int_{D_+} [\Delta_{(\xi_1,0)} f(x) - f_1'(x)\xi_1][x_1 - (x_1 \wedge lx_1x_2)]\mu(d\xi)$$

$$+ \int_{D_+} [\Delta_{(0,\xi_2)} f(x) - f_2'(x)\xi_2][lx_1x_2 - (x_1 \wedge lx_1x_2)]\mu(d\xi)$$

$$+ \int_{D_-} [\Delta_{(\xi_1,0)} f(x) - f_1'(x)\xi_1]x_1\mu(d\xi).$$

**7. Fluctuation limit theorems.** In this section we show that an affine process arises naturally from a limit theorem based on catalytic CBI-processes. By virtue of the characterizations given in the last section, we can establish the limit theorem in the sense of convergence in probability. Let $(\Omega, \mathscr{F}, \mathscr{F}_t, \mathbf{P})$ be a filtered probability space satisfying the usual hypotheses and let $B(\cdot)$, $N_0$ and $N_1$ be given as in the last section. Let $(a, (\alpha_{ij}), (\beta_{ij}), (b_j), m, \mu)$ be admissible parameters with $\beta_{22} < 0$ and $m_2(l_1) < \infty$. Let $\sigma_0 = \sqrt{a}$ and let $(\sigma_{ij})$ be a $(2 \times 2)$-matrix satisfying $(\alpha_{ij}) = (\sigma_{ij})(\sigma_{ij})^\tau$.

Let $\{\theta_k\}$ be a sequence such that $1 \leq \theta_k \to \infty$ as $k \to \infty$. For each $k \geq 1$ let $y_k(0)$ be an $\mathscr{F}_0$-measurable random variable and let $y_k(\cdot)$ be the solution of

(7.1)
$$y_k(t) = y_k(0) + \int_0^t (-\theta_k\beta_{22} + \beta_{21}x(s)\tilde{y}_k(s) + b_2\tilde{y}_k(s) + \beta_{22}y_k(s)) \, ds$$

$$+ \int_0^t \sigma_0\sqrt{2\tilde{y}_k(s)} \, dB_0(s) + \int_0^t \sigma_{21}\sqrt{2x(s)\tilde{y}_k(s)} \, dB_1(s)$$

$$+ \int_0^t \sigma_{22}\sqrt{2x(s)\tilde{y}_k(s)} \, dB_2(s) + \int_0^t \int_{D_+} \xi_2 \tilde{N}_0(ds, d\xi)$$

$$+ \int_0^t \int_0^{x(s-)\tilde{y}_k(s-)} \int_{D_+} \xi_2 \tilde{N}_1(ds, du, d\xi),$$

where $\tilde{y}_k(t) = y_k(t)/\theta_k$ and $x(\cdot)$ is defined by (6.6). When $\theta_k$ is sufficiently large, (7.1) is essentially a special form of (6.9). Then the pair $(x(\cdot), y_k(\cdot))$ is a catalytic CBI-process. Set $z_k(t) = y_k(t) - \theta_k$.

THEOREM 7.1. *Suppose that $z(0)$ is an $\mathscr{F}_0$-measurable random variable such that $\mathbf{E}[|z(0)|] < \infty$ and*

(7.2)
$$\lim_{k \to \infty} \mathbf{E}[|z_k(0) - z(0)|] = 0.$$

*Then $z_k(\cdot)$ converges in probability by the topology of $D([0,\infty),\mathbb{R})$ to the solution $z(\cdot)$ of*

$$z(t) = z(0) + \int_0^t (b_2 + \beta_{21}x(s) + \beta_{22}z(s))\,ds + \int_0^t \sqrt{2}\sigma_0\,dB_0(s)$$

(7.3)
$$+ \int_0^t \sigma_{21}\sqrt{2x(s)}\,dB_1(s) + \int_0^t \sigma_{22}\sqrt{2x(s)}\,dB_2(s)$$

$$+ \int_0^t \int_{D_+} \xi_2 \tilde{N}_0(ds,d\xi) + \int_0^t \int_0^{x(s-)} \int_{D_+} \xi_2 \tilde{N}_1(ds,du,d\xi).$$

PROOF. From (7.1) we get

$$\tilde{y}_k(t) = \tilde{y}_k(0) + \int_0^t (-\beta_{22} + \theta_k^{-1}\beta_{21}x(s)\tilde{y}_k(s) + \theta_k^{-1}b_2\tilde{y}_k(s) + \beta_{22}\tilde{y}_k(s))\,ds$$

$$+ \int_0^t \theta_k^{-1}\sigma_0\sqrt{2\tilde{y}_k(s)}\,dB_0(s) + \int_0^t \theta_k^{-1}\sigma_{21}\sqrt{2x(s)\tilde{y}_k(s)}\,dB_1(s)$$

(7.4)
$$+ \int_0^t \theta_k^{-1}\sigma_{22}\sqrt{2x(s)\tilde{y}_k(s)}\,dB_2(s) + \int_0^t \int_{D_+} \theta_k^{-1}\xi_2 \tilde{N}_0(ds,d\xi)$$

$$+ \int_0^t \int_0^{x(s-)\tilde{y}_k(s-)} \int_{D_+} \theta_k^{-1}\xi_2 \tilde{N}_1(ds,du,d\xi).$$

For $n \geq 1$ let $\tau_n = \inf\{s \geq 0 : x(s) \geq n\}$. Then $\tau_n \to \infty$ as $n \to \infty$. Under condition (7.2) we clearly have $\sup_{k\geq 1} \mathbf{E}[\tilde{y}_k(0)] < \infty$. By Proposition 5.2,

(7.5)
$$\sup_{k\geq 1} \mathbf{E}\left[\sup_{0\leq s\leq T} \tilde{y}_k(s\wedge\tau_n)\right] < \infty.$$

Let $\eta_{n,k}(t) = \tilde{y}_k(t\wedge\tau_n) - 1$. By (7.4) and Doob's martingale inequality we get

$$\mathbf{E}[|\eta_{n,k}(t)|] \leq \mathbf{E}[|\eta_{n,k}(0)|] + |\beta_{22}|\int_0^t \mathbf{E}[|\eta_{n,k}(s)|]\,ds$$

$$+ \theta_k^{-1}\int_0^t (b_2 + n|\beta_{21}|)\mathbf{E}[\tilde{y}_k(s\wedge\tau_n)]\,ds$$

$$+ \theta_k^{-1}(\sqrt{2}\sigma_0 + \sqrt{2n}\sigma_{21} + \sqrt{2n}\sigma_{22})\left(\int_0^t \mathbf{E}[\tilde{y}_k(s\wedge\tau_n)]\,ds\right)^{1/2}$$

$$+ \theta_k^{-1}\sqrt{tm_2(l_{12})} + 2\theta_k^{-1}tm_2(l_{12})$$

(7.6)
$$+ \theta_k^{-1}\sqrt{n\mu_2(l_{12})}\left(\int_0^t \mathbf{E}[\tilde{y}_k(s\wedge\tau_n)]\,ds\right)^{1/2}$$

$$+ 2n\theta_k^{-1}\mu_2(l_{12})\int_0^t \mathbf{E}[\tilde{y}_k(s\wedge\tau_n)]\,ds,$$

where $m_2$ and $\mu_2$ denote, respectively, the projections of $m$ and $\mu$ to $\mathbb{R}$. An application of Gronwall's inequality shows that

$$(7.7) \qquad \mathbf{E}[|\eta_{n,k}(t)|] = \mathbf{E}[|\tilde{y}_k(t \wedge \tau_n) - 1|] \to 0$$

as $k \to \infty$. From (7.1) we see that $z_k(\cdot)$ satisfies

$$(7.8) \begin{aligned} z_k(t) = {}& z_k(0) + \int_0^t (b_2 \tilde{y}_k(s) + \beta_{21} x(s) \tilde{y}_k(s) + \beta_{22} z_k(s))\, ds \\ & + \int_0^t \sigma_0 \sqrt{2\tilde{y}_k(s)}\, dB_0(s) + \int_0^t \sigma_{21} \sqrt{2x(s)\tilde{y}_k(s)}\, dB_1(s) \\ & + \int_0^t \sigma_{22} \sqrt{2x(s)\tilde{y}_k(s)}\, dB_2(s) + \int_0^t \int_{D_+} \xi_2 \tilde{N}_0(ds, d\xi) \\ & + \int_0^t \int_0^{x(s-)\tilde{y}_k(s-)} \int_{D_+} \xi_2 \tilde{N}_1(ds, du, d\xi). \end{aligned}$$

Let $\zeta_{n,k}(t) = z_k(t \wedge \tau_n) - z(t \wedge \tau_n)$. Then we have

$$\begin{aligned} \zeta_{n,k}(t) = {}& \zeta_{n,k}(0) + \int_0^{t \wedge \tau_n} (b_2 + \beta_{21} x(s)) \eta_{n,k}(s)\, ds \\ & + \beta_{22} \int_0^{t \wedge \tau_n} \zeta_{n,k}(s)\, ds + \int_0^{t \wedge \tau_n} \sqrt{2}\sigma_0 (\sqrt{\tilde{y}_k(s)} - 1)\, dB_0(s) \\ & + \int_0^{t \wedge \tau_n} \sigma_{21} \sqrt{2x(s)}(\sqrt{\tilde{y}_k(s)} - 1)\, dB_1(s) \\ & + \int_0^{t \wedge \tau_n} \sigma_{22} \sqrt{2x(s)}(\sqrt{\tilde{y}_k(s)} - 1)\, dB_2(s) \\ & + \int_0^{t \wedge \tau_n} \int_{x(s-)}^{x(s-)\tilde{y}_k(s-)} \int_{D_+} \xi_2 \tilde{N}_1(ds, du, d\xi). \end{aligned}$$

By Itô's formula,

$$\begin{aligned} e^{-\beta_{22} t} \zeta_{n,k}(t) = {}& \zeta_{n,k}(0) + \int_0^{t \wedge \tau_n} e^{-\beta_{22} s}(b_2 + \beta_{21} x(s)) \eta_{n,k}(s)\, ds \\ & + \int_0^{t \wedge \tau_n} \sqrt{2}\sigma_0 e^{-\beta_{22} s}(\sqrt{\tilde{y}_k(s)} - 1)\, dB_0(s) \\ & + \int_0^{t \wedge \tau_n} \sigma_{21} e^{-\beta_{22} s} \sqrt{2x(s)}(\sqrt{\tilde{y}_k(s)} - 1)\, dB_1(s) \\ & + \int_0^{t \wedge \tau_n} \sigma_{22} e^{-\beta_{22} s} \sqrt{2x(s)}(\sqrt{\tilde{y}_k(s)} - 1)\, dB_2(s) \\ & + \int_0^{t \wedge \tau_n} \int_{x(s-)}^{x(s-)\tilde{y}_k(s-)} \int_{D_+} e^{-\beta_{22} s} \xi_2 \mathbb{1}_{\{|\xi_2| \le 1\}} \tilde{N}_1(ds, du, d\xi) \end{aligned}$$

$$+ \int_0^{t\wedge\tau_n}\int_{x(s-)}^{x(s-)\tilde{y}_k(s-)}\int_{D_+} e^{-\beta_{22}s}\xi_2\mathbb{1}_{\{|\xi_2|>1\}}N_1(ds,du,d\xi)$$

$$- \int_0^{t\wedge\tau_n} e^{-\beta_{22}s}x(s)[\tilde{y}_k(s)-1]\,ds\int_{D_+}\xi_2\mathbb{1}_{\{|\xi_2|>1\}}\mu(d\xi).$$

By Doob's inequality we get

$$\mathbf{E}\left[\sup_{0\leq s\leq t}|e^{-\beta_{22}s}\zeta_{n,k}(s)|\right]$$

$$\leq \mathbf{E}[\zeta_{n,k}(0)] + \int_0^t e^{-\beta_{22}s}(b_2+n|\beta_{21}|)\mathbf{E}[|\eta_{n,k}(s)|]\,ds$$

$$+ 4\sqrt{2}\sigma_0\left(\int_0^t e^{-2\beta_{22}s}\mathbf{E}[(\sqrt{\tilde{y}_k(s\wedge\tau_n)}-1)^2]\,ds\right)^{1/2}$$

$$+ 4\sqrt{2n}(\sigma_{22}+\sigma_{21})\left(\int_0^t e^{-2\beta_{22}s}\mathbf{E}[(\sqrt{\tilde{y}_k(s\wedge\tau_n)}-1)^2]\,ds\right)^{1/2}$$

$$+ 4\sqrt{n\mu_2(l_{12})}\left(\int_0^t e^{-2\beta_{22}s}\mathbf{E}[|\eta_{n,k}(s)|]\,ds\right)^{1/2}$$

$$+ 2n\mu_2(l_{12})\int_0^t e^{-\beta_{22}s}\mathbf{E}[|\eta_{n,k}(s)|]\,ds,$$

where

$$\mathbf{E}[(\sqrt{\tilde{y}_k(s\wedge\tau_n)}-1)^2]\leq \mathbf{E}[|\tilde{y}_k(s\wedge\tau_n)-1|]=\mathbf{E}[|\eta_{n,k}(s)|].$$

Then (7.2) and (7.7) imply that

(7.9) $$\mathbf{E}\left[\sup_{0\leq s\leq t}|e^{-\beta_{22}s}\zeta_{n,k}(s)|\right]\to 0$$

as $k\to\infty$. For any $\varepsilon>0$, $\eta>0$ and $T\geq 0$ we first choose $n$ so that $\mathbf{P}\{\tau_n\leq T\}\leq \varepsilon/2$. In view of (7.9), there is some $k_0$ so that

$$\mathbf{P}\left\{\sup_{0\leq s\leq T}|\zeta_{n,k}(s)|\geq \eta\right\}\leq \eta^{-1}\mathbf{E}\left[\sup_{0\leq s\leq T}|\zeta_{n,k}(s)|\right]\leq \varepsilon/2$$

for every $k\geq k_0$. It then follows that

$$\mathbf{P}\left\{\sup_{0\leq s\leq T}|z_k(s)-z(s)|\geq \eta\right\}\leq \mathbf{P}\{\tau_n\leq T\}+\mathbf{P}\left\{\sup_{0\leq s\leq T}|\zeta_{n,k}(s)|\geq \eta\right\}\leq \varepsilon.$$

Then $z_k(\cdot)$ converges to $z(\cdot)$ in probability by the topology of $D([0,\infty),\mathbb{R})$.
□

Clearly, the pair $(x(\cdot),z(\cdot))$ defined by (6.6) and (7.3) is an affine process with nonnegative jumps. In other words, Theorem 7.1 gives an interpretation of a particular class of affine processes in terms of catalytic CBI-processes. To consider general affine processes, we assume the following decompositions of the parameters:

$$
\text{(7.10)} \quad \begin{aligned} \sigma_0 &= \sigma_0^+ - \sigma_0^-, & \sigma_{2j} &= \sigma_{2j}^+ - \sigma_{2j}^-, \\ b_2 &= b_2^+ - b_2^-, & \beta_{21} &= \beta_{21}^+ - \beta_{21}^-. \end{aligned}
$$

Let $x(\cdot)$ be defined by (6.6) and let $y_k^\pm(\cdot)$ be the solutions of the equations

$$
\begin{aligned}
y_k^+(t) = {}& y_k^+(0) + \int_0^t (-\theta_k \beta_{22} + \beta_{21}^+ x(s) \tilde{y}_k^+(s) + b_2^+ \tilde{y}_k^+(s) + \beta_{22} y_k^+(s))\, ds \\
& + \int_0^t \sigma_0^+ \sqrt{2\tilde{y}_k^+(s)}\, dB_0(s) + \int_0^t \sigma_{21}^+ \sqrt{2x(s)\tilde{y}_k^+(s)}\, dB_1(s) \\
& + \int_0^t \sigma_{22}^+ \sqrt{2x(s)\tilde{y}_k^+(s)}\, dB_2(s) + \int_0^t \int_{D_+} \xi_2 \tilde{N}_0(ds, d\xi) \\
& + \int_0^t \int_0^{x(s-)\tilde{y}_k^+(s-)} \int_{D_+} \xi_2 \tilde{N}_1(ds, du, d\xi),
\end{aligned}
\tag{7.11}
$$

$$
\begin{aligned}
y_k^-(t) = {}& y_k^-(0) + \int_0^t (-\theta_k \beta_{22} + \beta_{21}^- x(s) \tilde{y}_k^-(s) + b_2^- \tilde{y}_k^-(s) + \beta_{22} y_k^-(s))\, ds \\
& + \int_0^t \sigma_0^- \sqrt{2\tilde{y}_k^-(s)}\, dB_0(s) + \int_0^t \sigma_{21}^- \sqrt{2x(s)\tilde{y}_k^-(s)}\, dB_1(s) \\
& + \int_0^t \sigma_{22}^- \sqrt{2x(s)\tilde{y}_k^-(s)}\, dB_2(s) - \int_0^t \int_{D_-} \xi_2 \tilde{N}_0(ds, d\xi) \\
& - \int_0^t \int_0^{x(s-)\tilde{y}_k^-(s-)} \int_{D_-} \xi_2 \tilde{N}_1(ds, du, d\xi),
\end{aligned}
\tag{7.12}
$$

where $\tilde{y}_k^\pm(t) = y_k^\pm(t)/\theta_k$. We may regard $(x(\cdot), y_k^+(\cdot), y_k^-(\cdot))$ as a *catalytic CBI-process with a pair of reactant processes*. Set $z_k^\pm(t) = y_k^\pm(t) - \theta_k$ and $z_k(t) = z_k^+(t) - z_k^-(t) = y_k^+(t) - y_k^-(t)$.

THEOREM 7.2. *Suppose that $z^+(0)$ and $z^-(0)$ are $\mathscr{F}_0$-measurable random variables such that $\mathbf{E}[|z^+(0)| + |z^-(0)|] < \infty$ and*

$$
\text{(7.13)} \quad \lim_{k \to \infty} \mathbf{E}[|z_k^+(0) - z^+(0)| + |z_k^-(0) - z^-(0)|] = 0.
$$

*Then $z_k(\cdot)$ converges in probability by the topology of $D([0,\infty), \mathbb{R})$ to the solution $z(\cdot)$ of*

$$z(t) = z(0) + \int_0^t (b_2 + \beta_{21} x(s) + \beta_{22} z(s))\,ds + \int_0^t \sqrt{2}\sigma_0\,dB_0(s)$$

(7.14)
$$+ \int_0^t \sigma_{21}\sqrt{2x(s)}\,dB_1(s) + \int_0^t \sigma_{22}\sqrt{2x(s)}\,dB_2(s)$$
$$+ \int_0^t \int_D \xi_2 \tilde{N}_0(ds, d\xi) + \int_0^t \int_0^{x(s-)} \int_D \xi_2 \tilde{N}_1(ds, du, d\xi),$$

where $z(0) = z^+(0) - z^-(0)$.

PROOF. By Theorem 7.1, the sequence $(z_k^+(\cdot), z_k^-(\cdot))$ converges in probability by the topology of $D([0,\infty), \mathbb{R}^2)$ to the solution $(z^+(\cdot), z^-(\cdot))$ of

$$z^+(t) = z^+(0) + \int_0^t (b_2^+ + \beta_{21}^+ x(s) + \beta_{22} z^+(s))\,ds + \int_0^t \sqrt{2}\sigma_0^+\,dB_0(s)$$
$$+ \int_0^t \sigma_{21}^+ \sqrt{2x(s)}\,dB_1(s) + \int_0^t \sigma_{22}^+ \sqrt{2x(s)}\,dB_2(s)$$
$$+ \int_0^t \int_{D_+} \xi_2 \tilde{N}_0(ds, d\xi) + \int_0^t \int_0^{x(s-)} \int_{D_+} \xi_2 \tilde{N}_1(ds, du, d\xi),$$
$$z^-(t) = z^-(0) + \int_0^t (b_2^- + \beta_{21}^- x(s) + \beta_{22} z^-(s))\,ds + \int_0^t \sqrt{2}\sigma_0^-\,dB_0(s)$$
$$+ \int_0^t \sigma_{21}^- \sqrt{2x(s)}\,dB_1(s) + \int_0^t \sigma_{22}^- \sqrt{2x(s)}\,dB_2(s)$$
$$- \int_0^t \int_{D_-} \xi_2 \tilde{N}_0(ds, d\xi) - \int_0^t \int_0^{x(s-)} \int_{D_-} \xi_2 \tilde{N}_1(ds, du, d\xi).$$

It is simple to check that $z(\cdot) = z^+(\cdot) - z^-(\cdot)$ solves (6.7). That proves the theorem. □

The pair $(x(\cdot), z(\cdot))$ defined by (6.6) and (7.14) is an affine process with admissible parameters $(\sigma_0^2, (\alpha_{ij}), (\beta_{ij}), (b_j), m, \mu)$. Then the above theorem establishes a connection between catalytic CBI-processes and affine processes. This result is of interest since the studies of catalytic branching processes and affine processes have been undergoing rapid developments in recent years with rather different motivations; see, for example, [6] and [10].

**Acknowledgments.** We would like to thank Professors T. Shiga and S. Watanabe for helpful discussions and Professor Z. Q. Chen for useful comments on the literature. We are grateful to the two referees for pointing out a number of typos.

AFFINE MARKOV PROCESSES 41[23] Rhyzhov, Y. M. and Skorokhod, A. V. (1970). Homogeneous branching processes with a finite number of types and continuous varying mass. *Theory Probab. Appl.* **15** 704–707. MR0288860

[24] Sato, K. (1999). *Lévy Processes and Infinitely Divisible Distributions.* Cambridge Univ. Press. MR1739520

[25] Sharpe, M. J. (1988). *General Theory of Markov Processes.* Academic Press, New York. MR0958914

[26] Schmuland, B. and Sun, W. (2001). On the equation $\mu_{t+s} = \mu_s * T_s\mu_t$. *Statist. Probab. Lett.* **52** 183–188. MR1841407

[27] Shiga, T. and Watanabe, S. (1973). Bessel diffusions as a one-parameter family of diffusion processes. *Z. Wahrsch. Verw. Gebiete* **27** 37–46. MR0368192

[28] Watanabe, S. (1969). On two dimensional Markov processes with branching property. *Trans. Amer. Math. Soc.* **136** 447–466. MR0234531

[29] Yamada, T. and Watanabe, S. (1971). On the uniqueness of solutions of stochastic differential equations. *J. Math. Kyoto Univ.* **11** 155–167. MR0278420
School of Mathematics and Statistics  
Carleton University  
1125 Colonel By Drive  
Ottawa  
Canada K1S 5B6  
E-mail: ddawson@math.carleton.ca  
URL: lrsp.carleton.ca/directors/dawson/

School of Mathematical Sciences  
Beijing Normal University  
Beijing 100875  
People's Republic of China  
E-mail: lizh@email.bnu.edu.cn  
URL: math.bnu.edu.cn/~lizh/